\documentclass[12pt]{amsart}
\usepackage{amsthm,amsfonts, amssymb, amscd}

\newcommand{\C}{{\mathbb C}}

\newcommand{\F}{{\mathcal F}}
\newcommand{\T}{{\mathcal T}}
\newcommand{\E}{{\mathcal E}}

\newcommand{\m}{{\mathfrak m}}
\newcommand{\g}{{\mathfrak g}}
\newcommand{\z}{{\mathfrak z}}

\DeclareMathOperator{\Spec}{Spec}
\DeclareMathOperator{\Td}{Td}
\DeclareMathOperator{\ch}{ch}
\DeclareMathOperator{\GL}{GL}

\newtheorem{thm}{Theorem}
\newtheorem*{intro-thm}{Theorem}

\newtheorem{prop}[thm]{Proposition}

\newtheorem{lemma}[thm]{Lemma}

\newtheorem{lemma-definition}[thm]{Lemma-Definition}

\newtheorem{cor}[thm]{Corollary}
\theoremstyle{remark}
\newtheorem{remark}[thm]{Remark}
\newtheorem{definition}[thm]{Definition}
\begin{document}
\numberwithin{thm}{section}
\title[Nonabelian localization]{Nonabelian localization in equivariant $K$-theory 
and Riemann-Roch for quotients}
\author{Dan Edidin}
\address{
Department of Mathematics,
University of Missouri,
Columbia, MO 65211
}
\email{edidin@math.missouri.edu}
\author{William Graham}
\address{
Department of Mathematics,
University of Georgia,
Boyd Graduate Studies Research Center,
Athens, GA 30602
}
\email{wag@math.uga.edu}
\thanks{The first author was supported by N.S.A. and the 
second by N.S.F}
\maketitle
\begin{abstract}
We prove a localization formula in equivariant algebraic $K$-theory for
an arbitrary complex algebraic group acting with finite stabilizer
on a smooth algebraic space. This extends to non-diagonalizable
groups the localization formulas of H.A. Nielsen \cite{Nie:74} and 
R. Thomason \cite{Tho:92}

As an application we give a Riemann-Roch formula for quotients of
smooth algebraic spaces by proper group actions. This formula extends
previous work of B. Toen \cite{Toe:99} and the authors \cite{EdGr:03}.
\end{abstract}
\tableofcontents
\section{Introduction}
Equivariant $K$-theory was developed in the late 1960's by Atiyah and
Segal as a tool for the proof of the index theorem for elliptic
operators invariant under the action of a compact Lie group. In the
late 1980's and early 1990's Thomason constructed an algebraic
equivariant $K$-theory modeled on Quillen's earlier construction of
higher $K$-theory for schemes.   

In both the topological and algebraic contexts equivariant $K$-theory
is studied using its structure as a module for the representation ring
of the group $G$.  The fundamental theorem of equivariant $K$-theory
is the {\it localization theorem} for actions of diagonalizable
groups.  We describe a version of this theorem in the complex algebraic
setting. Let $R(G)$ be the representation ring of $G$ tensored with
$\C$. If $X$ is a $G$-space we let $G(G,X)$ be the equivariant 
$K$-groups of the category of
$G$-equivariant coherent sheaves, also tensored with $\C$.  If $G$ is
diagonalizable and $h \in G$,
let $\iota: X^h \rightarrow X$ be the inclusion of the fixed locus of
$h$.  Let $\m_h \subset R(G)$ denote the maximal ideal of
representations whose virtual characters vanish at $h$.  The
localization theorem states that the natural map
$$
\iota_*: G(G,X^h)_{\m_h} \rightarrow G(G,X)_{\m_h},
$$
is an isomorphism.  
Much of the power of this theorem comes from the
fact that if $X$ is regular, then so is $X^h$, and the localization
isomorphism has an explicit inverse, arising from the
self-intersection formula for the regular embedding $X^h
\stackrel{\iota} \to X$: If $\alpha \in G(G,X)_{\m_h}$ is an element
of localized equivariant $K$-theory then
\begin{equation} \label{eqn.awesome}
\alpha = \iota_*\left(  \frac{ \iota^*\alpha}{  \lambda_{-1}(N_\iota^*) } \right).
\end{equation}
Here $N_\iota^*$ is the conormal bundle to
$\iota$, and $\lambda_{-1}(N_\iota^*)$ is defined
to be the element in $K$-theory corresponding to the formal sum 
$\sum_{l = 0}^{\text{rank} N_\iota^*}(-1)^{l} \Lambda^l(N_\iota^*).$

The formula of Equation \eqref{eqn.awesome} is extremely useful
because it reduces global calculations to those on the fixed locus.
It has been applied in a wide range of contexts. For example, the
localization theorem on the flag variety $G/B$ can be used to give a
proof of the Weyl character formula.  In \cite{EdGr:03} we used the
localization theorem to prove a Kawasaki-Riemann-Roch formula for
quotients by diagonalizable group actions (similar ideas had been
introduced earlier by Atiyah \cite{Ati:74}).

If we try to generalize Equation \eqref{eqn.awesome} to the
nonabelian case we immediately run into the problem that $X^h$ is not
in general $G$-invariant.  However, the locus $X_{\Psi} = \overline{G X^h}$ is
$G$-invariant; it is the closure of the union of the fixed point loci of
elements in the conjugacy class $\Psi$ of $h$.  
Let $\m_{\Psi} \subset
R(G)$ denote the maximal ideal of representations whose virtual
characters vanish on $\Psi$, and let $i :X_{\Psi} \rightarrow X$
denote the inclusion.  Then $i_*\colon G(G,X_{\Psi})_{\m_{\Psi}} \to
G(G,X)_{\m_{\Psi}}$ is an isomorphism (see Theorem \ref{thm.weakloc};
this is a variant of a result of Thomason \cite{Tho:92}, adapting a
result of Segal \cite{Seg:68b} from topological $K$-theory).
Unfortunately, as Thomason observed \cite{Tho:88}, $X_{\Psi}$ can be
singular even when $X$ is smooth, so the self-intersection formula
does not apply.  To obtain a nonabelian version of Equation
\eqref{eqn.awesome} new ideas are needed.

Although $X^h$ is not $G$-invariant, it is $Z$-invariant, where $Z =
{\mathcal Z}_G(h)$ is the centralizer in $G$ of $h$. In \cite{VeVi:02}
Vezzosi and Vistoli proved that if $G$ acts on $X$ with finite
stabilizers then there is an isomorphism between a localization of
$G(Z,X^h)$ and a localization of $G(G,X)$. (Note that their theorem
holds in arbitrary characteristic.)  In this paper we work with the
added hypothesis that the projection $f$ from the global stabilizer $S_X =
\{(g,x) | gx = x\}$ to $X$ is a finite morphism. Our main
result states that there is a natural pushforward $\iota_! \colon
G(Z,X^h) \to G(G,X)$, such that when $X$ is smooth, the following
formula holds for $\alpha \in G(G,X)_{\m_\Psi}$:
\begin{equation} \label{eqn.irock}
\alpha = \iota_! \left( \frac{ \lambda_{-1}(({\mathfrak g}/{\mathfrak
z})^*)\cap (\iota^!\alpha)_h }{ \lambda_{-1}(N_\iota^*) } \right).
\end{equation}
Here $\iota^!$ is the composition of the restriction map $G(G,X) \to G(Z,X)$
with the pullback $G(Z,X) \stackrel{\iota^*} \to G(Z,X^h)$; 
$(\iota^!\alpha)_h$ is the image of 
$\iota^! \alpha$ in $G(Z,X^h)_{\m_h}$, and $\mathfrak{g}$,
$\mathfrak{z}$ are the Lie algebras of $G$ and $Z$ respectively.

To prove this result, we use an equivalent formulation involving the
global stabilizer. Let $S_{\Psi} \subset S_X$ be the closed subspace
of pairs $(g,x)$ with $g \in \Psi$.  The finite map $f:
S_{\Psi} \rightarrow X$ has image $X_{\Psi}$, but unlike $X_{\Psi}$,
the space $S_{\Psi}$ is regular (if $X$ is).  There is a natural
identification of $G(G,S_{\Psi})$ with $G(Z,X^h)$, and the map
$\iota_!$ is identified with the pushforward $f_*$ in $G$-equivariant
$K$-theory.  Moreover, the natural map $f: S_{\Psi} \rightarrow X$ is
a local complete intersection morphism, so it has a normal bundle
$N_f$.  There is a distinguished "central summand"
$G(G,S_{\Psi})_{c_{\Psi}}$ of $G(G,S_{\Psi})$; if $\beta \in
G(G,S_{\Psi})$, we let $\beta_{c_{\Psi}}$ denote the component of
$\beta$ in the central summand.  Equation \eqref{eqn.irock} is
equivalent to the statement that if  $\alpha \in G(G,X)_{\m_{\Psi}}$, then
\begin{equation} \label{eqn.irock2}
\alpha = f_* \left( \frac{(f^*
\alpha)_{c_{\Psi}}}{\lambda_{-1}(N_f^*)} \right).
\end{equation}
This formula looks similar to the formula that would hold if $f$ were
a regular embedding (the only change would be to replace $(f^*
\alpha)_{c_{\Psi}}$ by $f^*{\alpha}$).  However, that formula is not
correct, and indeed, the main difficulty in proving \eqref{eqn.irock2} is
that $f$ is not a regular embedding, so we cannot apply the
self-intersection formula.  The proof given here is less direct; we
first prove the result when $G$ is a product of general linear groups,
and then use a change of groups argument to deduce the general case.

The main application of Equations \eqref{eqn.irock} and
\eqref{eqn.irock2} is to give refined formulas for the Todd classes of
sheaves of invariant sections on quotients of smooth algebraic spaces.
If $X$ is a smooth, separated, algebraic space and $G$ is an algebraic
group acting properly (and thus with finite stabilizer) then the
theorem of Keel and Mori \cite{KeMo:97} implies that there is a
(possibly singular) geometric quotient $Y = X/G$. For such quotients
there is a map in $K$-theory $\pi_G\colon G(G,X) \to G(Y)$ induced by
the exact functor which takes a $G$-equivariant coherent sheaf to its
subsheaf of invariant sections (Lemma \ref{lem.invariants}).  Let
$\tau_Y \colon G_0(G,X) \to A_*(Y)$ be the Riemann-Roch map defined by
Baum, Fulton, MacPherson \cite{BFM:75,Ful:84}.  If $\alpha \in
G_0(G,X)$ then we obtain explicit expressions (Theorems
\ref{t.rrtheorem1} and \ref{t.rrtheorem2}) for $\tau_Y(
\pi_G(\alpha_{\Psi}))$ in terms of the restriction of $\alpha$ to a
class in $G_0(Z,X^h)$ where $h \in \Psi$ is any element. If we sum
over all conjugacy classes $\Psi$ we obtain formulas for
$\tau_Y(\pi_G(\alpha))$.  When the quotient is quasi-projective our formulas
for $\tau_Y(\pi_G(\alpha))$ can be deduced
from the Riemann-Roch formula for stacks due to B. Toen
\cite{Toe:99}. Our method of proof is quite different, and makes no
essential use of stacks.

The proof of our Riemann-Roch theorem is essentially the same as the
proof for diagonalizable $G$ given in \cite{EdGr:03}, with the
nonabelian localization theorem of this paper in place of the
localization theorem for diagonalizable groups.  A key element of the
proof in \cite{EdGr:03} was the fact that, if $G$ is diagonalizable,
$X^h$ is $G$-invariant and we may define
an $h$-action on
$G_0(G,X^h)$ which we called ``twisting by $h$''.  
Intuitively, this twist
comes from the $h$-action on the sections of any 
$G$-equivariant coherent sheaf on $X^h$.  
In the nonabelian setting one can still twist 
by a central element, so there is an action of $h$
on $G_0(Z, X^h)$.  This can be viewed as a twist of $G(G,S_{\Psi})$,
which intuitively comes from the tautological action of the element
$g$ on the fiber at $(g,x)$ of any $G$-equivariant vector bundle on
$S_{\Psi}$.  Versions of this twist and the global stabilizer appear
in the Riemann-Roch theorems of Kawasaki and Toen, and motivated our
approach to the localization theorem and Riemann-Roch theorem.
Interestingly, in the Riemann-Roch formula obtained from
\eqref{eqn.irock2}, one might expect a term involving $(f^*
\alpha)_{c_{\Psi}}$.  However, we prove that 
the contributions from $(f^* \alpha)_{c_{\Psi}}$ and $f^*
\alpha$ are equal, so our formula does not mention
the central summand.

In this paper we work over $\C$ and tensor all $K$-groups with
$\C$. The reason we do this is so that we can identify the
representation ring $R(G)$ with class functions on $G$; this idea,
which goes back to Atiyah and Segal, allows us to directly relate
$G$-equivariant $K$-theory to conjugacy classes in $G$.  By working
over $\C$ we hope that the geometric techniques used to prove our main
results are not obscured by technical details.  

Nevertheless, we believe that
versions of the nonabelian localization and Riemann-Roch theorems
should hold over an arbitrary algebraically closed field provided we
assume that all stabilizer groups are reduced. In this situation,
instead of localizing at maximal ideals $\m_h \in R(G) \otimes \C$
where $h \in G$ has finite order, we may localize at the
multiplicatively closed set $S_H$ defined on p.~10 of \cite{VeVi:02},
where $H$ is the cyclic group generated by $h$.  In a different direction,
there should be topological versions of these results for actions of
compact Lie groups.  This will be pursued
elsewhere.

We would like to thank Michel Brion for pointing out an error in an earlier
version of the paper. We are also grateful to the referee
for several helpful suggestions.

\subsection{Conventions and notation}  \label{s.conventions}

We work entirely over the ground field of complex numbers
$\C$. All algebraic spaces are assumed to be 
of finite type over $\C$. For a reference on the theory
of algebraic spaces, see the book \cite{Knu:71}.  All 
algebraic groups are assumed to be linear. A basic reference
for the theory of algebraic groups is Borel's book
\cite{Bor:91}.

If $G$ is an algebraic group then ${\mathcal Z}(G)$ denotes
the center of $G$. If $h \in G$ is any element then
${\mathcal Z}_G(h)$ denotes the centralizer of $h$ in $G$.
The conjugacy class of $h$ in $G$ is denoted $C_G(h)$.
The map $G \to C_G(h)$, $g \mapsto ghg^{-1}$ identifies 
$C_G(h)$ with the homogeneous space $G/{\mathcal Z}_G(h)$.

If $G$ is an algebraic group then $R(G)$ denotes
the representation ring of $G$ tensored with $\C$.

\subsubsection{Group actions} \label{ss.group}
Let $G$ be an algebraic group acting on an algebraic space $X$.
We consider three related conditions on group actions.

(i) We say that $G$ acts {\it properly} if the map
$$G\times X \to X \times X, \; (g,x) \mapsto (x,gx)$$ is
proper.

Let $G \times X \to X \times X$ be the map we just defined
and let $X \to X \times X$ be the diagonal. Then
$S_X = G \times X \times_{X \times X} X$ is called the global
stabilizer. As a set, $S_X = \{(g,x)| gx =x\}$.

(ii) We say that $G$ acts with {\it finite stabilizer} if 
the projection $S_X \to X$ is a finite morphism.

(iii) We say that $G$ acts with {\it finite stabilizers}
if the projection $S_X \to X$ is quasi-finite; i.e. 
for every point $x \in X$ the isotropy group $G_x$ is finite.

Since $G$ is affine, the map $G \times X \to X \times X$
is finite if it is proper. It follows that any proper
action has finite stabilizer. If $G$ acts
properly on $X$, then the diagonal morphism of $X$ 
factors as the composition of two proper maps
$X \stackrel{(e,1_X)} \to G \times X \to X \times X$.
This implies that the diagonal is a closed embedding,
so $X$ is automatically separated.

\section{Groups, representation rings, and conjugacy classes} \label{s.repring}

This section collects a number of facts about algebraic groups,
representation rings, and conjugacy classes which are used in the
proof of the localization theorem.  We have included proofs of some
results that are essentially known but are difficult to find in the
literature for groups that are not connected or semisimple.

\subsection{The representation ring and class functions}
Let $G$ be an algebraic group over $\C$. The group $G$ is called
reductive if the radical of its identity component $G_0$ is a torus.
Any reductive algebraic group is the complexification of a maximal
compact subgroup \cite[Theorem 8, p.~244]{OnVi:90}.  (This result can be deduced
from the corresponding result for connected groups, using facts about
maximal compact subgroups of Lie groups with finitely many connected
components \cite[Theorem XV.3.1]{Hoc:65}.)  By the unipotent
radical of $G$ we mean the unipotent radical of the identity component
of $G$; this is a normal subgroup of $G$, and $G$ modulo its unipotent
radical is reductive.  Any complex algebraic group $G$ has a Levi
subgroup $L$.  This means $G$ is the semidirect product of $L$ and the
unipotent radical of $G$ (\cite[Chapter 6, Theorem 4]{OnVi:90}); $L$
is necessarily reductive.

For a general algebraic group $G$, let $\hat{G}$ denote the set of isomorphism
classes of irreducible (finite-dimensional) algebraic representations of $G$.  Let
$\C[G]$ denote the coordinate ring of $G$, and $\C[G]^G$ the
ring of class functions, i.e., the functions on $G$ which are invariant under
conjugation.  There is a map $R(G) \rightarrow \C[G]^G$ which
takes $[V]$ to $\chi_V$, where $V$ is a representation of $G$ and
$\chi_V$ its character.  

Let $V$ be a 
representation of $G$,
and let $V^*$ denote the dual representation.  There is a map
\begin{equation} \label{e.repfunction}
V^* \otimes V \rightarrow \C[G]
\end{equation}
taking $\lambda \otimes v$ to the function $\lambda/v$ defined by
$$
\lambda/v(x) = \lambda(xv)
$$
for $x \in G$.  This map is $G \times G$-equivariant, where $G \times G$ acts
on $V^* \otimes V$ via the $G$-action on each factor, and $G \times G$ acts
on functions by
$$
( (g_1,g_2) \cdot f)(x) = f(g_1^{-1} x g_2),
$$
for $g_1,g_2,x \in G$ and $f \in \C[G]$.  Functions of the form $\lambda / v$
are called representative functions; we denote the algebra they generate by
$\T(G)$.  If $G$ is a linear algebraic group, then the action of $G \times G$ on
$\C[G]$ is locally finite, and this can be used to show that $\T(G) = \C[G]$.  If $G$ is the complexification of $K$, then the map $\T(G) \rightarrow \T(K)$ is an isomorphism (cf. \cite{BrtD:95}).  As a consequence, we
obtain the algebraic Peter-Weyl theorem, attributed to Hochschild and Mostow:

\begin{prop}
If $G$ is a complex reductive algebraic group, then the map
\begin{equation} \label{e.repfunction2}
\oplus_{V \in \hat{G}}V^* \otimes V \rightarrow \C[G]
\end{equation}
is an isomorphism as representations of $G \times G$. 
\end{prop}

This follows from the usual Peter-Weyl theorem for compact groups
(\cite[p.~201]{Ros:02}),
in view of the isomorphism $\C[G] \rightarrow \T(K)$.

In this paper, we frequently use the following result.  For connected
groups, it can be proved using restriction to a maximal torus.

\begin{prop} \label{prop.repring}
If $G$ is reductive, the map $R(G) \rightarrow \C[G]^G$
taking a representation to its character is an isomorphism. In particular, for any
$G$,
$R(G)$ is a finitely generated algebra over $\C$.
\end{prop}

\begin{proof}
View $G$ as embedded diagonally in $G \times G$.  The isomorphism
\eqref{e.repfunction2} induces an isomorphism of $G$-invariants in the
source and target.  If $V$ is any representation of $G$, then $(V^*
\otimes V)^G = \mbox{Hom}_G(V,V)$, and if $V$ is irreducible
representation, Schur's lemma implies that this is $1$-dimensional,
spanned by the identity map $\mbox{id}_V$.  But $\phi(\mbox{id}_V) =
\chi_V$, proving the first statement of the proposition.  The second statement
follows because the representation rings of $G$ and any Levi factor are isomorphic.
\end{proof}

As another application of the relationship between $G$ and $K$, we have the following proposition.

\begin{prop} \label{prop.finite}
If $H \hookrightarrow G$ is an embedding of groups, then $R(H)$ is a finite $R(G)$-module.
\end{prop}

\begin{proof}
If $H_1$ is a Levi subgroup of $H$ then the restriction map $R(H) \rightarrow R(H_1)$ is an isomorphism,
so if necessary replacing $H$ by $H_1$, we may assume $H$ is reductive.  
Let $G'$ be the quotient of $G$ by its unipotent radical; then the natural map $R(G') \rightarrow R(G)$ is an isomorphism.  Moreover, the kernel of $H \rightarrow G'$ is
a unipotent normal subgroup of $H$, hence is trivial.  Hence the map $H \rightarrow G'$ is injective.
Therefore, if necessary replacing $G$ by $G'$, we may assume $G$ is reductive.

Let $L$ be a maximal compact subgroup of $H$; we may assume $L$ is contained
in a maximal compact subgroup $K$ of $G$.  Then the restriction maps
$R(G) \rightarrow R(K)$ and $R(H) \rightarrow R(L)$ are isomorphisms
(cf. \cite{BrtD:95}).  Therefore the proposition
follows from the analogous result of \cite{Seg:68a}
for compact groups.
\end{proof}

\subsection{Conjugacy classes and representation rings}
If $H$ is a subgroup of $G$, let $C_H(g)$ denote the $H$-conjugates of
$g$, i.e., the image of the map $H \rightarrow G$ taking $h$ to
$hgh^{-1}$.

\begin{prop} \label{prop.conjclosed}
Let $G$ be a reductive algebraic group.  A conjugacy class $\Psi$ in $G$ is closed if and only
if it is semisimple.
\end{prop}

\begin{proof}
In any algebraic group, the conjugacy class of a semisimple element is
closed \cite[Theorem~9.2]{Bor:91}.  Conversely, suppose $\Psi$ is a
closed conjugacy class.  Let $g = su$ be the Jordan decomposition of
some $g \in \Psi$.  We want to show that $u = 1$, or in other words,
that $C_G(s) = C_G(su)$.  A general result of Mumford about reductive
groups acting on affine schemes (\cite[Ch.~1, Cor.~1.2]{MFK:94})
implies that given two closed disjoint conjugation-invariant subsets
of $G$, there is a class function which is $1$ on one subset and $0$
on the other.  Therefore, to show that $C_G(s) = C_G(su)$ it suffices
to show that any class function takes the same value on $C_G(s)$ as on
$C_G(su)$.  But this holds because it is true for characters, which
span the space of class functions.
\end{proof}

Note that if $G$ is not reductive, there may be closed conjugacy
classes which are not semisimple, for example if $G = {\mathbb G}_a$.

If $\Psi = C_G(g)$ is a semisimple conjugacy class, let $\m_\Psi
\subset R(G)$ be the ideal of virtual characters which vanish on
$\Psi$.  Then $\m_\Psi$ is the kernel of the homomorphism of
$\C$-algebras $R(G) \to \C$ defined by the property that $[V] \mapsto
\chi_V(\Psi)$. Since the trivial character does not vanish on $\Psi$,
$\m_\Psi$ is a maximal ideal of $R(G)$.

\begin{prop} \label{prop.specarbitrary}
Let $G$ be an arbitrary algebraic group. 
The assignment $\Psi \mapsto \m_\Psi$ gives
a bijection between the set of semisimple conjugacy
classes in $G$ and the maximal ideals in $R(G)$.
\end{prop}

\begin{proof}
First observe that the result holds if $G$ is reductive by Propositions
\ref{prop.repring} and \ref{prop.conjclosed}.
For general $G$, let $G = LU$ be a Levi decomposition and let
$r \colon R(G) \to R(L)$ be the restriction map. Since $U$ is
unipotent, $r$ is an isomorphism.

If $\m \subset R(G)$ is a maximal ideal then, since $L$ is reductive,
$r(\m)= \m_\Phi$ for some semisimple conjugacy class $\Phi = C_L(l)$
in $L$. If $\Psi = C_G(l)$ then $r(\m_\Psi) \subset \m_\Phi$. Thus
$\m_\Psi \subset \m$. But $\m_\Psi$ is maximal so it equals $\m$.

Suppose $\Psi_1$ and $\Psi_2$ are two semisimple conjugacy classes
in $G$ such that $\m_{\Psi_1} = \m_{\Psi_2}$. Since
the proposition holds for the reductive subgroup $L$, 
$r(\m_{\Psi_i}) = \m_\Phi$ for some
semisimple conjugacy class $\Phi \subset L$. It follows that
$\Psi_1 \cap L = \Psi_2 \cap L = \Phi$. This means that
$\Psi_1 \cap \Psi_2$ is nonempty, so $\Psi_1 = \Psi_2$.
\end{proof}

Given an embedding of groups $G \rightarrow H$, let $r\colon R(H) \to
R(G)$ be the restriction map.  Let $\Psi = C_H(h)$ be a semisimple
conjugacy class in $H$. By Proposition \ref{prop.finite}, the ideal $\m_\Psi
R(G)$ is contained in a finite number of maximal ideals of $\Spec
R(G)$. If $\Psi'$ is a conjugacy class in $\Psi \cap G$, then the
restriction of any virtual character vanishing on $\Psi$ also
vanishes on $\Psi'$; that is, $\m_\Psi R(G) \subset \m_{\Psi'}$. Thus,
$\Psi \cap G$ decomposes into a finite number of semisimple conjugacy
classes $\Psi_1, \ldots , \Psi_l$.

\begin{prop} \label{prop.conjdecomp}
Let $G \hookrightarrow H$ be an embedding of groups, and
let $\Psi'$ and $\Psi$ be semisimple conjugacy classes in $G$ and $H$, respectively.
Then
$\m_\Psi R(G) \subset \m_{\Psi'}$ if and only if
the conjugacy class $\Psi'$ is contained in $(\Psi \cap G)$.
\end{prop}
\begin{proof} 
The if direction follows from the discussion 
immediately preceding the statement of the proposition.
Conversely, if $\Psi' $ is a conjugacy class in $G$
not contained in $\Psi \cap G$, then $\Psi'$ is disjoint from the conjugacy
classes $\Psi_1, \ldots , \Psi_l$ in $\Psi \cap G$.  By Proposition
\ref{prop.specarbitrary}, $\m_{{\Psi}'}$ is distinct from each of the
$\m_{\Psi_i}$.  Therefore
there is a virtual
character $f \in \m_{\Psi'}$ such $f$ is not in some $\m_{\Psi_i}$.
Thus, the restriction of $f$ to $\Psi \cap G$ is not zero, so $f \notin \m_\Psi R(G)$.
\end{proof}

\begin{remark} \label{rem.pickbill}
Proposition \ref{prop.conjdecomp} implies that if $\Psi \subset H$
is a semisimple conjugacy class then $R(G)_{\m_{\Psi}}$ is
a semilocal ring with maximal ideals $\m_{\Psi_{1}}R(G)_{\m_{\Psi}},
\ldots \m_{\Psi_{l}}R(G)_{\m_{\Psi}}$. 
\end{remark}

As noted above, if $G$ is a subgroup of $H$ and $g \in
G$, the intersection of $C_H(g)$ with $G$ may consist of more than one
$G$-conjugacy class in $G$.  The following result shows that it is
possible to find embeddings where this does not occur.

\begin{prop}\label{prop.goodembedding}
Suppose $G$ is an algebraic group and $\Psi = C_G(g)$ is a 
semisimple conjugacy class in $G$.   There is an embedding $G \rightarrow H$,
where $H = \prod_i \GL_{n_i}$, such that $C_H(g) \cap G = C_G(g)$.
\end{prop}

\begin{proof} Since $R(G) = R(L)$, where $L$ is a Levi factor of $G$,
$R(G)$ is Noetherian. 
Therefore we can find a finite set $f_1,\ldots f_l$ of elements which 
generate the maximal ideal $\m_\Psi \subset R(G)$.
Each function $f_i$ can be written as a finite sum $f_i = \sum_j 
a_{ij}\chi_{ij}$ where $\chi_{ij}$ is the character of a $G$-module 
$V_{ij}$. Let $V$ be a faithful representation of $G$ and
set $H = \GL(V) \times \prod_{i,j} \GL(V_{ij})$. Then
$G$ embeds as a subgroup of $H$ since it embeds as a subgroup of
the first factor, $\GL(V)$. 

Let $g'$ be an element of $G$ such that $g'$ is conjugate to $g$
in $H$. Since $H$ is a direct product, the image of $g$ and $g'$ 
in each $\GL(V_{ij})$ must have the same trace. Thus, $\chi_{ij}(g) = 
\chi_{ij}(g')$ for all $i,j$. Hence $f_{i}(g)= f_i(g') = 0$ for each generator
$f_i$ of $\m_\Psi$ (since $g$ is, by definition, an element
of the conjugacy class $\Psi$).  
Therefore, by Proposition \ref{prop.specarbitrary},
$g'$ is in the conjugacy class of $g$ in $G$.
\end{proof}

\section{Equivariant $K$-theory} \label{s.ekt}
This section contains some $K$-theoretic results needed for
the proof of our main result, the nonabelian localization theorem.

\subsection{Basic facts and notation}
Let $G$ be an algebraic group acting on an algebraic space $X$. We use
the notation ${\tt coh}^G_X$ to denote the abelian category of
$G$-equivariant coherent ${\mathcal O}_X$-modules.  Let
$$G(G,X) = \oplus_{i = 0}^\infty G_i(G,X) \otimes \C$$ where
$G_i(G,X)$ is the $i$-th Quillen $K$-group of the category ${\tt
coh}^G_X$.  Since all our coefficients are taken to be complex, we
will simply write $G_0(G,X)$ (rather than $G_0(G,X) \otimes \C$) for
the Grothendieck group of $G$-equivariant coherent sheaves, tensored
with $\C$.  Likewise, we write $K_0(G,X)$ for the Grothendieck ring of
$G$-equivariant bundles, also tensored with $\C$. When $X$ is a smooth
scheme, Thomason's equivariant resolution theorem implies that
$K_0(G,X) = G_0(G,X)$ (this is true even without tensoring with $\C$).
We use analogous notation in the non-equivariant setting, writing
$\tt{coh}_X$ for the category of coherent ${\mathcal O}_X$-modules,
and write $G(X)$, $G_0(X)$, and $K_0(X)$ for the non-equivariant
versions of $K$-theory.

If ${\mathcal E}$ is a $G$-equivariant locally free sheaf
then the assignment ${\mathcal F} \mapsto {\mathcal E \otimes F}$
defines an exact functor ${\tt coh}_X^G \to {\tt coh}_X^G$. 
This implies that there is an action of $K_0(G,X)$ on $G(G,X)$.

If $X$ and $Y$ are $G$-spaces and $p \colon X \to Y$ is a $G$-map, then
there is a pullback $p^* \colon K_0(G,X) \to K_0(G,X)$. When $Y =
\Spec \C$, this pullback makes $G(G,X)$ an $R(G) = K_0(G,\Spec \C)$
module. 

If $p\colon X \to Y$ is a $G$-equivariant proper morphism then there
is a pushforward $p_* \colon G(G,X) \to G(G,Y)$ which is an 
$R(G)$-module homomorphism \cite[1.11-12]{Tho:86i}.
The pushforward is defined on the level 
of Grothendieck groups by the formula
$p_*[{\mathcal F}] = 
\sum (-1)^{i}[R^i p_* {\mathcal F}].$

If $p \colon X \to Y$ 
is flat and $G$-equivariant 
then there is a pullback $p^*\colon G(G,X) \to G(G,Y)$
which is also an $R(G)$-module homomorphism.
It is induced by the exact functor which takes a coherent
sheaf ${\mathcal F}$ to its pullback $p^*{\mathcal F}$.

More generally, if $X$ is a regular algebraic space then Vezzosi and Vistoli 
proved \cite[Theorem A4]{VeVi:02} that $G(G,X)$ is isomorphic
to the Waldhausen $K$-theory of the category ${\mathcal W}_{3,X}$
of complexes of flat quasi-coherent $G$-equivariant ${\mathcal O}_X$
modules with bounded coherent cohomology. It follows
that if $p \colon X \to Y$ is a map of regular algebraic
spaces then there is a pullback $p^* \colon G(G,X) \to G(G,Y)$ which
is an $R(G)$-module homomorphism.
If $p$ is a regular embedding of smooth algebraic spaces
then there is a self-intersection formula for $ \alpha \in G(G,X)$:
$$p^* p_* \alpha = \lambda_{-1}(N_p^*) \cap \alpha,$$
where $N_p^*$ is the conormal bundle to map $p$ and 
$\lambda_{-1}(N_p^*) \in K_0(G,X)$ is the formal
sum $\sum_{i = 0}^{rk\; N_p^*} (-1)^i
[\Lambda^i(N_p^*)]$.
This fact is proved in the course of the proof of Theorem 3.7
of \cite{VeVi:02}.

\subsection{Morita equivalence} \label{s.morita}
In this paper we make extensive use of a particular instance of Morita equivalence, which
we briefly describe.  
If $Z \subset G$ is a closed subgroup and $X$ is a $Z$-space
then we may consider the $G \times Z$-space $G \times X$
where $(k,z) \cdot (g,x) = (kgz^{-1},zx)$. 
We write $G \times_Z X$
for the quotient of $G \times X$ by the free action
of the subgroup $1 \times Z$. The space
$G \times_Z X$ will often be referred to as a {\it mixed space}. 
If the $Z$ action on $X$
is the restriction of a $G$ action, then
the automorphism of $G \times X$ given by $(g,x) \mapsto (g,gx)$
induces an isomorphism of quotients $G\times_Z X \to G/Z \times X$.

Since the actions of $G$ and $Z$ on $G \times X$ commute, the action
of $G \times 1$ on $G \times X$ descends to an action on the quotient
$G \times_Z X$.  The Morita equivalence we use is the equivalence of categories
between the category of $Z$-equivariant coherent sheaves on $X$ and
the category $G$-equivariant coherent sheaves on $G \times_Z X$.  The
equivalence is given by pulling a $Z$-module on $X$ back to $G \times
X$ to obtain a $G \times Z$-module on $G \times X$ and then taking the
subsheaf of $1 \times Z$-invariant sections to obtain a $G$-module on $G
\times_Z X$. 

\begin{remark} \label{rem.tantomixedbund} 
When $X$ is a point, Morita equivalence between the categories of
$G$-equivariant coherent sheaves on $G/Z$ and $Z$-modules is obtained
by taking the fiber of a sheaf at the identity coset $Z$.  Under this
equivalence the tangent bundle of $G/Z$ corresponds to the $Z$-module
${\mathfrak g}/{\mathfrak z}$, where ${\mathfrak g}$ and ${\mathfrak
z}$ denote the Lie algebras of $G$ and $Z$ respectively 
\cite[Proposition 6.7]{Bor:91}. If $X$ is a
$G$-space then the identification $G \times_Z X = G/Z \times X$
identifies the relative tangent bundle of the projection $p \colon G \times_Z X
\to X$ with the bundle $G \times_Z (X \times {\mathfrak g}/{\mathfrak
z})$.
\end{remark}

\begin{remark} \label{rem.moritaequiv}
The Morita equivalence of categories above induces an
$R(G)$-module isomorphism in $K$-theory
$G(Z,X) \to G(G, G \times_Z X)$. Here
$R(G)$ acts on $G(Z,X)$ via the restriction map
$R(G) \to R(Z)$. This observation will be used repeatedly
in the sequel.
\end{remark}

\subsection{Localization in equivariant $K$-theory} \label{sec.loc}
In this section we extend the localization theorem of
\cite[Theorem 2.2]{Tho:92} to arbitrary algebraic groups.
The proof is similar to Thomason's.  However, because
we work with ideals in $R(G) \otimes \C$ which may have zero
intersection with the integral representation ring, we cannot
directly quote 
his results.

\begin{thm} \label{thm.weakloc} Let $G$ be an
algebraic group acting on an algebraic space $X$. Let $\Psi = C_G(h)$
be a semisimple conjugacy class and let $X_\Psi$ be the closure
of $GX^h$ in $X$.

(a) The proper
pushforward
$$i_*\colon  G(G,X_\Psi) \to G(G,X)$$ 
is an isomorphism of $R(G)$-modules after localizing at $\m_\Psi$.

(b) If $X$ is smooth and $h \in {\mathcal Z}(G)$ (so $\Psi = h$ and $X^h$ is
$G$-invariant)
then the map of $R(G)$-modules 
$$\cap\;\lambda_{-1}(N_i^*) \colon G(G,X^h) \to G(G,X^h)$$ is invertible after localizing
at $\m_h$, and if $\alpha \in G(G,X)_{\m_h}$, then
\begin{equation} \label{eqn.prettygood}
\alpha = i_*\left(\lambda_{-1}(N_i^*)^{-1} \cap i^*\alpha\right).
\end{equation}
Here the notation $\lambda_{-1}(N_i^*)^{-1} \cap i^*\alpha$
means the image of $i^*\alpha$ under the inverse of
the isomorphism $\cap \; \lambda_{-1}(N_i^*)$.
\end{thm}
\begin{remark} If $X$ is a smooth scheme then $X^h$ is as well.
By Thomason's equivariant resolution theorem we
may identify 
$K_0(G,X^h)$ with $G_0(G,X^h)$ \cite[Theorem 5.7]{Tho:87}
and view $(\lambda_{-1}(N_i^*))^{-1}$
as an element in $G(T,X^h)_{\m_h}$.
\end{remark}

\begin{proof}
{\em Step 1: $G = T$ is a torus.}
If $T$ is a torus then $C_G(h) = h$. Using Noetherian induction
and the localization long exact sequence 
\cite[Theorem 2.7]{Tho:87} it suffices
to prove that $G(T,X)_{\m_h} = 0$ whenever $X^h$ is empty.
Since $T$ is diagonalizable, we can, by Proposition \ref{prop.repring},
identify $R(T)$ with the coordinate
ring of $T$, and $\m_h$ with the maximal ideal of $h \in T$.
If $T' \subset T$ is a closed subgroup then $R(T')_{\m_h} = 0$ 
unless $h \in T'$. By Thomason's generic slice theorem 
\cite[Theorem 4.10]{Tho:86d}
there is a $T$-invariant open set $U$, a closed subgroup $T'$ acting
trivially on $U$ and a $T$-equivariant isomorphism
$U \simeq T/T' \times U/T$. Since
$T'$ acts trivially on $U$, and $X^h$ is empty, $h \notin T'$.
By Morita equivalence, 
$G(T,U) = G(T',U/T)$. But $T'$ acts trivially
on $U/T$, so $G(T',U/T) = R(T') \otimes G(U/T)$. Thus
$G(T',U/T)_{\m_h} = 0$.  Applying
Noetherian induction and using the localization long exact sequence
we conclude that $G(T,X)_{\m_h} =0$. This proves part (a).

Next we must show that when $X$ is smooth, 
the multiplication
map 
$$\cap\;\lambda_{-1}(N_i^*) \colon G(T,X^h)_{\m_h} \to G(T,X^h)_{\m_h}$$ is an
isomorphism.  This may be done using Noetherian induction on $X^h$.  Again
we apply 
Thomason's generic slice theorem. Thus we may assume that $X^h = T/T'
\times X^h/T$. As above
$G(T,X^h)_{\m_h} = R(T')_{\m_h} \otimes G(X^h/T)$.

Let ${\mathcal N}_{x}$ be the fiber of $N_i^*$ over a point $x\in
X^h$. Since $T'$ acts trivially on $X^h$, 
${\mathcal N}_x$ is a $T'$-module.  As
in the proof of \cite[Lemma~3.2]{Tho:92}, the identification
$G(T,X^h)_{\m_h} = R(T')_{\m_h} \otimes G(X^h/T)$ implies that the action
of $\lambda_{-1}(N_i^*)$ on $G(T,X^h)_{\m_h}$ is invertible if and
only if $\lambda_{-1}({\mathcal N}_x)$ is invertible in $R(T')_{\m_h}$
for some $x$ in each connected component of $X^h$; 
i.e. $\lambda_{-1}({\mathcal N}_x) \notin \m_h$.
Now ${\mathcal N}_x$ decomposes into $1$-dimensional eigenspaces for
the action of $T'$ so we can write $\lambda_{-1}({\mathcal N}_x)=
\prod_{l=1}^s (1 - \chi_l)$, where $\chi_1, \ldots , \chi_s$ are (not
necessarily distinct) characters of $T'$.

Let $H$ be the closure of the cyclic subgroup of $T$ generated by
$h$. Then $X^H = X^h$.  Since $X^h = T/T' \times X^h/T$, $T'$ is the
biggest subgroup of $T$ acting trivially on $X^h = X^H$. Thus, we see
that $H \subset T'$ and the characters in the $T'$-module
decomposition restrict to characters of $H$.  None of these characters
can be trivial on $H$ since the normal space to $X^H$ at $x$ is the
quotient of the tangent space $T_{x,X}$ by the invariant subspace
$T_{x,X}^H$. Since the cyclic subgroup generated by $h$ is dense in
$H$, we see that, for each $\chi_l$, we have $\chi_l(h^n) =
\chi_l(h)^n \neq 1$ for some exponent $n$.  Thus $(1 - \chi_l) \notin
\m_h$.  Therefore, the action of $\lambda_{-1}(N_i^*)$ on
$G(T,X^h)_{\m_h}$ is invertible.

Finally, the formula in Equation \eqref{eqn.prettygood}
can be deduced as follows.
Since $i_*$ is surjective, $\alpha = i_* \beta$ for some
$\beta \in G(T,X^h)_{\m_h}$. Thus,
\begin{eqnarray*}
i^*\alpha &  = & i^* i_* \beta\\
 & = & \lambda_{-1}(N_i^*)\cap\beta
\end{eqnarray*}
where the second equality follows from the self-intersection
formula in equivariant $K$-theory.
Since the action of 
$\lambda_{-1}(N_i^*)$ is invertible after localizing at $\m_h$,
the formula follows.\\

{\em Step 2.  $G$ is connected and reductive.}   We use a standard
reduction to a maximal torus argument.  Let $T$ be any maximal torus
in $G$ and let $B$ be a Borel 
subgroup of $G$ containing $T$. Since $B/T$ is isomorphic to affine space,
the restriction $G(B,X) \to G(T,X)$ is an isomorphism 
\cite[Proof of Theorem 1.13]{Tho:88}.  
The same proof implies that the
projection $p \colon G \times_B X \to X$ induces a pullback $p^!\colon
G(G,X) \to G(G,G \times_B X)= G(T,X)$ and a pushforward $p_!\colon
G(T,X) \to G(G,X)$, with the properties that $p_! 1 =1$, and
if $\beta \in K_0(G,X)$, then $p_! (p^!\beta \cap \alpha) = \beta \cap
p_!\alpha$.  In particular, $p^!$ is a split monomorphism. Moreover,
standard functorial properties of equivariant $K$-theory imply that
$p^!$ and $p_!$ are functorial for $G$-equivariant morphisms.

To prove (a), as in the torus case it suffices to prove that if $X^h$ is empty
then $G(G,X)_{\m_{\Psi}} = 0$.  Since $G(G,X)$ embeds in $G(T,X)$ it
actually suffices to prove that $G(T,X)_{\m_{\Psi}} = 0$.  Let
$h_1, h_2, \ldots h_n$ be the conjugates of $h$ contained in the
maximal torus $T$.   By Remark \ref{rem.pickbill},
$R(T)_{\m_{\Psi}}$ is a semilocal ring whose maximal ideals
are the $\m_{h_i} R(T)_{\m_{\Psi}}$.  Hence, 
if $M$ is any $R(T)$-module, to show that
$M_{\m_{\Psi}} = 0$ it suffices to show that $M_{\m_{h_i}} = 0$ for all $i$.

Since each $h_i$ is $G$-conjugate to $h$, each
$X^{h_i}$ is empty.   By the torus case, $G(T,X)_{\m_{h_i}} =
0$ for each $h_i$.  This implies that $G(G,X)_{\m_{\Psi}} = 0$, proving (a).

Now suppose that $h \in {\mathcal Z}(G)$. To avoid confusion,
we write $\m_h^G$ for the maximal ideal in $R(G)$ corresponding
to the one-element conjugacy class $h \in G$ and $\m_h^T$ the maximal 
ideal in $R(T)$ corresponding to the one-element conjugacy class $h \in T$.
By Remark \ref{rem.pickbill}, $R(T)_{\m_h^G} = R(T)_{\m_h^T}.$ If $X$ is smooth
then, by Step 1, the action of $\lambda_{-1}(N_i^*)\in K_0(G,X^h)$ on 
$G(T,X^h)_{\m_h^G} = G(T,X^h)_{\m_h^T}$ is invertible. 
Thus the action of $\lambda_{-1}(N_i^*)$
on $G(G,X^h)_{\m_h^G} = p_!(G(T,X^h)_{\m_h^G})$ is as well.

Once we know that the action of $\lambda_{-1}(N_i^*)$ is invertible
after localizing at $\m_h^G$, the formula of 
equation \eqref{eqn.prettygood} follows from the self-intersection formula.\\

{\em Step 3: $G$ is arbitrary.}
As in the previous steps, to prove part (a) it suffices to show that if
$\Psi = C_G(h)$
and $X^h$ is empty, then $G(G,X)_{\m_\Psi} = 0$.  By Proposition
\ref{prop.goodembedding}, there is
an embedding of $G$ into a product of general linear groups $Q$ such
that if $\Psi_1 = C_Q(h)$, then $\Psi_1 \cap Q = \Psi$.  This implies that
if $X^h$ is empty then so is $(Q \times_G X)^h$. Thus, by 
Remark \ref{rem.moritaequiv} and Step 2, we conclude
that $G(G,X)_{\m_{\Psi_1}} = 0$. Since $\m_{\Psi_1} R(G) \subset \m_\Psi$
it follows that $G(G,X)_{\m_\Psi}= 0$ as well.

Now suppose that $h$ is central in $G$; then $h$ is central in $Q$ as well.
Let $i \colon X^h \to X$ be the inclusion of the fixed locus
and let $\iota_h \colon Q \times_G X^h \to Q \times_G X^h$
be the corresponding inclusion of mixed spaces. By Morita
equivalence and Remark \ref{rem.pickbill}, it suffices to show that the action
of $\lambda_{-1}(N_{\iota}^*)$ on $G(Q, Q \times_G X^h)$ is invertible
after localizing at $\m_h^Q \subset R(Q)$. This follows from Step 2 and the following
lemma.

\begin{lemma} \label{lem.changefixed}
Let $G$ be a closed subgroup of an algebraic
group $Q$ and let $h \in G \cap {\mathcal Z}(Q)$. If
$X$ is a smooth $G$-space then 
$(Q \times_G X^h) = (Q \times_G X)^h$ 
as closed subspaces of $Q \times_G X$.
\end{lemma}
\begin{proof}[Proof of Lemma \ref{lem.changefixed}]
It is clear that $Q \times_G X^h \subset (Q \times_G X)^h$, so
we need only show the reverse inclusion. 
Since $Q \times_G X^h$ and $(Q \times_G X)^h$ are closed
smooth subspaces of the algebraic space $Q \times_G X$
it suffices to show that they have the same
closed points (since we work over the algebraically closed field $\C$).

A point corresponding to the $G$-orbit of $(q,x)\in Q \times X$
is fixed by $h$ if and only if $h(q,x) = (hq,x)$ is in the same $G$
orbit as $(q,x)$. This means that there is an element $g \in G$
such that $(hq,x) = (q g^{-1},g x)$. Thus $g$ fixes $x$ and 
$g^{-1} = q^{-1}h q$. Since $h \in {\mathcal Z}(Q)$ this means $g= h^{-1}$.
Since $h$ and $h^{-1}$ have the same fixed locus we conclude that
$x \in X^{h}$; i.e. $(q,x) \in Q \times X^h$.
\end{proof}
This concludes the proof of Theorem \ref{thm.weakloc}.
\end{proof}

\subsection{Decomposition of equivariant $K$-theory} \label{sec.decomp}
Let $G$ be an algebraic group acting on an algebraic space
$X$. Assume that $G$ acts with finite stabilizers. In
this case, there is a decomposition of $G(G,X)$ into
a direct sum of pieces, which we now describe.

Since $X$ is assumed to be Noetherian there is a finite
set of conjugacy classes $\Phi_1, \ldots \Phi_m$ of elements
of finite order such that $X^g$ is nonempty if and only
if $g \in \Phi_i$ for some $i$ \cite[Theorem 5.4]{VeVi:02}.
\begin{prop} \label{prop.decomp}
With the assumptions above, the localization maps
$G(G,X) \to G(G,X)_{\m_{\Phi_i}}$ induce a direct sum decomposition
$$G(G,X) = \oplus_i G(G,X)_{\m_{\Phi_i}}.$$
\end{prop}
\begin{proof}
By \cite[Remark 5.1]{EdGr:00}, there is a ideal $J \subset R(G)$ that
annihilates $G(G,X)$ and such that $R(G)/J$ is supported at a finite
number of points\footnote{If $G(G,X)$ is a finitely
generated $R(G)$-module this also follows from Theorem
\ref{thm.weakloc} and \cite[Theorem 5.4]{VeVi:02}.} of $\Spec R(G)$.  This implies that
$G(G,X) = \oplus_i G(G,X)_{\m_{\Phi_i}}$ for some set of semi-simple
conjugacy classes $\{\Phi_1, \Phi_2, \ldots \Phi_m\}$.  If $X^h$ is
empty for $h \in \Phi$, then by Theorem \ref{thm.weakloc},
$G(G,X)_{\m_\Psi} = 0$.
\end{proof}

\begin{remark} Suppose that $G$ acts on $X$ with finite stabilizers.
If $\alpha \in G(G,X)$, then we denote the component of 
$\alpha$ in the summand $G(G,X)_{\m_\Psi}$ by $\alpha_\Psi$. 
Note that if $\beta \in K_0(G,X)$ then 
$$(\beta \cap \alpha)_\Psi = \beta \cap (\alpha_\Psi).$$
Also, suppose that $f$ is a $G$-equivariant morphism of algebraic
spaces such that $G$ acts with finite stabilizers on the source
and target.  
If $f$ is proper morphism, then $f_*(\alpha_\Psi) = (f_* \alpha)_\Psi$.
Likewise, if $f$ is flat or a map of regular 
algebraic spaces then $f^*\alpha_\Psi = (f^* \alpha)_\Psi$. 
These basic facts follow immediately from the fact
that $f_*$ and $f^*$ are $R(G)$-module homomorphisms.
They will be used repeatedly in the proof of Theorem \ref{thm.localization}.
\end{remark}

Let $X$ be a $Z$-space where $Z$ acts with finite stabilizers and let $Z
\subset G$ be an embedding of $Z$ into another algebraic
group $G$. Morita equivalence (Section \ref{s.morita}) identifies
$G(G, G \times_Z X)$ with $G(Z,X)$ giving an
$R(Z)$-module structure on $G(G,G \times_Z X)$.  As a result
we may obtain a more refined decomposition of $G(G, G \times_Z X)$.
\begin{prop} \label{prop.changedecomp}
Let $Z$ act on $X$ with finite stabilizers.  
If $\Psi$ is a semisimple conjugacy class in $G$ and
$(\Psi \cap Z)$ decomposes into the union of conjugacy classes
$\Psi_1, \ldots , \Psi_l$, then 
$$G(G,G\times_Z X)_{\m_\Psi} = \oplus_{i=1}^l G(G,G \times_Z X)_{\m_{\Psi_i}}.$$
\end{prop}
\begin{proof}
By Remark \ref{rem.pickbill}, $R(Z)_{\m_\Psi}$ is a semi-local ring
with maximal ideals $\m_{\Psi_1},\ldots , \m_{\Psi_l}$ where
$\Psi_1, \ldots \Psi_l$ are the conjugacy classes in $\Psi \cap Z$.
The proposition follows.
\end{proof}

\subsection{Projection formulas for flag bundles}

In this section we prove some projection formulas for maps of flag
bundles which will be needed in the proof of the nonabelian
localization theorem. We begin with a lemma which is certainly known
in greater generality, but for which we lack a suitable reference.

\begin{lemma} \label{lem.projection} Let $p \colon P_1 \to P_2$ be a proper, flat, $G$-equivariant
map of quasi-projective schemes such that
action of $G$ is linearized with respect to an ample
line bundle on each $P_i$. Let $X$ be an algebraic space
with a $G$-action. Let $p_1 \colon P_1 \times X  \to P_1$
and $p_2 \colon P_2 \times X \to P_2$ be the projections. 
Set $\phi = (p \times 1) \colon P_1 \times X \to P_2 \times X$.
If $A \in K_0(G,P_1)$ and $\alpha \in G(G,P_2 \times X)$
then 
\begin{equation} \label{eqn.projection}
\phi_*(p_1^*A \cap \phi^*\alpha) = p_2^*p_*A \cap \alpha
\end{equation}
\end{lemma}

\begin{proof} We use
the ideas in the proof of the projection formula given
in \cite[Proposition 7.2]{Qui:73}.

By assumption the action
of $G$ on $P_1$ is linearized with respect to an ample
line bundle. Hence $K_0(G,P_1)$ is generated
by classes of equivariant vector bundles ${\mathcal E}$ with
$R^ip_*{\mathcal E} = 0$ for $i > 0$ \cite[Section 7.2]{Qui:73}. Thus
we may assume $A = [{\mathcal E}]$ and $R^ip_*{\mathcal E} = 0$ for $i>0$.
Since $p$ is flat it follows that $p_*{\mathcal E}$ is 
a locally free $G$-equivariant sheaf on $P_2$.
Thus, if ${\mathcal F}$ is any coherent sheaf on $P_2 \times X$
then $R^i\phi_*(p_1^*{\mathcal E} \otimes \phi^*{\mathcal F}) = 0$. (This can be
checked locally in the \'etale topology so we may assume $X$ is 
an affine scheme. The proof in that case is given on p. 59 of Srinivas's
book
\cite{Sri:96}.) Thus the functor 
$${\tt coh}^G_{(P_2 \times X)} \to {\tt coh}^G_{(P_2 \times
X)},\; \; {\mathcal F} \mapsto \phi_*(p_1^*{\mathcal E} \otimes 
\phi^*{\mathcal F})$$
is exact. This functor induces the endomorphism of $G(G,X)$
given by $\alpha \mapsto \phi_*(p_1^*A \cap \phi^*\alpha)$.

On the other hand, the endomorphism of $G(G,P_2 \times X)$
given by $\alpha \mapsto p_2^*p_*A \cap \alpha$ 
is induced by the exact functor
$${\tt coh}^G_{(P_2 \times X)} \to {\tt coh}^G_{(P_2 \times
X)},\; \; {\mathcal F} \mapsto p_2^*p_*{\mathcal E} \otimes {\mathcal F}.$$

Srinivas also proves that there is
a natural (and hence $G$-equivariant) isomorphism
$\phi_*(p_1^*{\mathcal E} \otimes \phi^*{\mathcal F}) \simeq 
\phi_*p_1^*{\mathcal E} \otimes {\mathcal F}$. Since 
$\phi_* p_1^*{\mathcal E}$ is naturally (and thus $G$-equivariantly)
isomorphic to $p_2^*p_*{\mathcal E}$ the two exact functors
we have defined are isomorphic. Therefore, the formula
of Equation \eqref{eqn.projection} holds.
\end{proof}

Let $G$ be a connected group 
and let $P \subset G$ be a parabolic subgroup containing a Borel
subgroup $B$ and having Levi factor $Z$. Choose a maximal torus $T
\subset Z$. Since $P$ is parabolic, $T$ is a maximal torus of $G$ as
well. Let $W(G,T) = N_G(T)/{\mathcal Z}_G(T)$ and $W(Z,T)= N_Z(T)/{\mathcal Z}_Z(T)$ be the Weyl groups of $G$ and $Z$
respectively.  If $X$ is a $G$-space
then we have projections
$$
G \times_B X\stackrel{p}{\rightarrow} G \times_P X\stackrel{q}{\rightarrow} X.
$$
Set
$\pi = q \circ p$.
The flag bundle projection formulas are given by the next result.

\begin{prop} \label{prop.tanpush}
If $\alpha \in G(G,X)$, $\beta \in G(G,G \times_P X)$ then the following
identities hold:\\

(i) $\pi_* (\lambda_{-1}(T^*_{\pi}) \cap \pi^*\alpha) = |W(G,T)|\alpha$.\\

(ii) $q_* (\lambda_{-1}(T^*_q) \cap q^*\alpha) = \frac{|W(G,T)|}
{|W(Z,T)|} \alpha$.\\

(iii) $p_* (\lambda_{-1}(T^*_{\pi}) \cap p^*\beta) = |W(Z,T)|
(\lambda_{-1}(T^*_q)\cap \beta)$.\\

\end{prop}

\begin{proof}
Since $G$ acts on $X$, the mixed spaces $G \times_B X$ and
$G \times_P X$ are isomorphic
to $G/B \times X$ and $G/P \times X$ respectively. The maps
$p,q, \pi$ are all smooth and projective and the 
bundles $T_{\pi}$, $T_{q}$ and $T_p$ are all obtained
by pullback from the smooth projective schemes $G/B$ or $G/P$.
Therefore, by Lemma \ref{lem.projection} we may assume that $X = \Spec \C$
and $\alpha = 1$.

Observe that since $B$ is a parabolic with Levi factor $T$, (i)
is a special case of (ii).  To prove (ii), it suffices to show that if $P \subset G$
is a parabolic with Levi factor $Z$ 
and $q \colon G/P \to \Spec \C$ is the projection,
then
$$q_*(\lambda_{-1}(T_{G/P}^*)) = |W(G,T)|/|W(Z,T)| \in R(G).$$
Since $R(G) \subset R(T)$ we may check this identity in $R(T)$. 
The torus  $T$ acts on $G/P$ with a finite number of fixed points
$p_1, \ldots , p_n$ where $n = |W(G,T)|/|W(Z,T)|$. Let 
$i_l \colon p_l \to G/P$ be the inclusion. The fixed points
are isolated and
$N_{i_l} = i_l^*T_{G/P}$.   We claim that
\begin{equation} \label{e.tanpush}
\lambda_{-1}(T_q^*) = \sum_{l=1}^n i_{l*} 1.
\end{equation}
Indeed, since $G/P$ has a $T$-equivariant cell decomposition,
\cite[Lemma 5.5.1]{ChGi:97} implies that $G_0(T,G/P) = K_0(T,G/P)$ 
is a free $R(T)$-module. 
In addition $R(T)$ is an integral domain, we may check \eqref{e.tanpush}
after localizing 
at any prime ideal in $R(T)$.   By Theorem \ref{thm.weakloc}, the equation holds
after localizing at $\m_a \in R(T)$ where $a\in T$ is any element
with $(G/P)^a = (G/P)^T$.  This proves \eqref{e.tanpush}.  Pushing forward
to $G(T,\Spec k) = R(T)$ completes the proof of (ii).

If $G$ is connected then $G$-equivariant $K$-theory is a summand in
$T$-equivariant $K$-theory \cite[Theorem 1.13]{Tho:88}.  Thus, to prove (iii) we may again work in $T$-equivariant $K$-theory.
The $T$ action on $G/B$ has $|W(G,T)|$ fixed points while
the $T$ action on $G/P$ has $|W(G,T)/W(Z,T)|$. For each
fixed point $P \in G/P$ the fiber
$p^{-1}(P)$ contains exactly $|W(Z,T)|$ of the $T$-fixed points in $G/B$.
Applying \eqref{e.tanpush} to $G/B$ and $G/P$, 
we see that $p_*\lambda_{-1}(T^*_{G/B}) = 
|W(Z,T)| \lambda_{-1}(T_{G/P}^*)$.
\end{proof}

\section{The global stabilizer and its equivariant $K$-theory} \label{s.globalstabilizer}

\subsection{General facts about the global stabilizer}
If $X$ is a $G$-space then $G$ acts on $G \times X$ by conjugation on
the first factor and by the original action on the second factor. The
global stabilizer $S_X \subset G \times X$ is a $G$-invariant subspace
and the projection $f \colon S_X \to X$ is $G$-equivariant.

\begin{definition} Let $\Psi$ be a semisimple conjugacy
class in $G$. Define $S_{\Psi} \subset S_X$ to be the inverse
image of $\Psi$ under the projection $S_X \to G$. 
Set theoretically, $S_\Psi = \{(g,x) | g \in \Psi
\text{ and } gx = x\}$. 
\end{definition}

\begin{remark}  \label{rem.finite}
Since semisimple conjugacy classes are closed,
$S_\Psi$ is closed in $S_X$. Thus, if $G$ acts with finite stabilizer
then the projection $f\colon S_\Psi \to X$ is also finite.

If $\Psi$ and $\Psi'$ are disjoint conjugacy classes then 
$S_\Psi \cap S_{\Psi^{'}} = \emptyset$.
When $G$ acts with finite stabilizers then 
$S_X$ is the disjoint sum of the closed subspaces
$S_{\Psi_{1}} \coprod S_{\Psi_{2}} \ldots \coprod S_{\Psi_{l}}$
where $\{\Psi_1, \ldots , \Psi_l\}$ is the set
of conjugacy classes whose elements have non-trivial stabilizer (these
conjugacy classes are necessarily semisimple because
in characteristic 0 any element of finite order is semisimple).
Hence $S_\Psi$ is also open in $S_X$ in this case.
\end{remark}

Let $h$ be an element of $\Psi$ and
let $Z = {\mathcal Z}_G(h)$ be the centralizer of $h$.
The map 
$G \times X^h \stackrel{\Phi_h} \to S_\Psi$ given by $(g,x) \mapsto (ghg^{-1},gx)$
is invariant under the free action of $Z$ on $G \times X^h$
given by $z(g,x) = (gz^{-1}, zx)$.
\begin{lemma} \label{lem.whatisspsi}
$G \times X^h \stackrel{\Phi_h} \to S_\Psi$ is a $Z$-torsor. Hence
$S_\Psi$ is identified with the quotient $G \times_Z X^h$.
\end{lemma}
\begin{proof} The map $G \to \Psi$, $g \mapsto ghg^{-1}$ identifies
$\Psi$ with $G/Z$. Thus, by base change, the map 
$\Phi_h \colon G \times X \to \Psi \times X$
given by $(g,x) \mapsto (ghg^{-1}, gx)$ is also a torsor. Since
$\Phi_h^{-1}(S_\Psi) = G\times X^h$ the lemma follows.
\end{proof}

\begin{remark} \label{rem.spsismooth}
If $X$ is smooth then 
the identification of $S_\Psi$ with $G \times_Z X^h$ implies that
$S_\Psi$ is smooth and 
the projection $f\colon S_\Psi \to X$ factors
as the composition of the regular embedding $S_\Psi \stackrel{i} \to G \times_Z X$
with the smooth projection morphism $G \times_Z X \stackrel{\eta}\to X$.
\end{remark}

If $X$ is smooth, 
the map $S_\Psi \stackrel{i} \to G \times_Z X$ is a regular
embedding but the projection $G \times_Z X \rightarrow Z$ is not proper.
This makes it difficult to use $G \times_Z X$ in the proof of our
main result.  However, suppose that $G$ is connected and $Z$ is a Levi factor of a parabolic subgroup $P$.  Let $\rho: G \times_Z X \rightarrow G \times_P
X$ be the projection, and let $j = \rho \circ i: S_\Psi \to G \times_P
X$.  The following result holds.

\begin{lemma} \label{lemma.zpsitoxp}
Let $G$ be connected and let $\Psi = C_G(h)$, where $h$ is semisimple.  Assume that
$Z = {\mathcal Z}_G(h)$ is a Levi factor of a parabolic subgroup $P$.
If $G$ acts on $X$ with finite stabilizer
then $j \colon S_{\Psi} \to G \times_P X$ is a regular embedding.
\end{lemma}
\begin{proof}
Since the morphisms of algebraic spaces 
$i$ and  $\rho$ are representable\footnote{We say that
a morphism of algebraic space $f\colon X \to Y$ is representable
if for any scheme $Y'$ and map $Y' \to Y$, the fiber product
$X \times_Y Y'$ is a scheme.} the composite $j = \rho \circ i$
is as well. By \cite[Proposition 19.1.1]{EGA4} any closed immersion
of regular schemes is a regular embedding. The property of being a
regular embedding is local for the \'etale topology
of the target so we may apply this proposition to representable morphisms
of algebraic spaces.
Since $S_\Psi$ and $G \times_P X$ are smooth we
are reduced to showing that $j$ is a closed immersion.

The finite map $f\colon S_{\Psi} \to X$ factors as $p \circ j$
where $p\colon G \times_P X \to X$ is the projection,
so $j$ is finite and thus closed. To prove that it is an
immersion
we will show that is unramified and injective on geometric points.
These two conditions suffice
because by \cite[Cor.~18.4.7]{EGA4} an unramified morphism
(Zariski) locally factors as the composition of an \'etale morphism
with a closed immersion. If, in addition, the map is injective
on geometric points
then the \'etale morphism is an open immersion by \cite[Theorem 17.9.1]{EGA4}.
Hence these conditions will imply that our morphism is locally the composition of an open immersion
with a closed immersion.

To show that $j$ is unramified it suffices to show
that the map $f = j \circ q \colon S_\Psi \to X$ is unramified.
Since $S_\Psi$ is open in $S_X$ we need only show that
the projection $S_X \to X$ is unramified. This last
fact follows from the fact the fiber of $S_X \to X$ over a geometric
point $x \to X$ is a group scheme over $x$. Since we work in characteristic
0 this must be reduced.

Next we show that $j$ is injective on geometric points.  Consider the
morphism $G \times X^h \to G \times_Z X^h \to G \times_P X$.  Suppose
that $(g_1,x_1)$ and $(g_2,x_2)$ have the same image in $G \times_P
X$. By definition this means that there is an element $p \in P$ such
that $x_2 = px_1$ and $g_2 = g_1 p^{-1}$.  Since $P$ is a semidirect
product of $Z$ and $U$ with $U$ unipotent and normal, we can write $p
= uz$ with $u \in U$ and $z \in Z$.  Thus $x_2 = u z x_1$; we want to
show that $u = 1$.  Because $z x_1$ and $x_2$ are both in $X^h$, it
suffices to show that if $u$ is an element of $U$ with $x$ and $ux$ both
in $X^h$, then $u = 1$.  This last assertion follows because the
assumptions imply that $h u x = h u h^{-1} x = u x$, or $u^{-1}
(huh^{-1})$ fixes $x$.  Now, $u^{-1} (huh^{-1})$ is unipotent since it
is in $U$.  However, this element has nontrivial fixed locus, so by
Remark \ref{rem.finite} it is semisimple.  Hence $u^{-1} (huh^{-1}) =
1$, so $u \in Z$.  Since $Z \cap U = \{ 1 \}$ we obtain $u = 1$, as
desired.
\end{proof}

\subsection{Localization theorems related to the global stabilizer}
Let $\Psi = C_G(h)$ be a semisimple conjugacy class, and let $Z =
{\mathcal Z}_G(h)$.  Recall that by Lemma \ref{lem.whatisspsi},
$\Phi_h: G \times X^h \rightarrow S_{\Psi}$ is a $Z$-torsor.

By Morita equivalence, $G(G,S_\Psi)$ is isomorphic to
$G(Z,X^h)$, so 
$G(G,S_\Psi)$ is an $R(Z)$-module.  The restriction homomorphism $R(G) \rightarrow
R(Z)$ is compatible with the $R(G)$ and $R(Z)$-module structures on
$G(G,S_\Psi)$ (cf. Section \ref{s.morita}).
Since $h$ is central in $Z$
it is a one-element conjugacy class with corresponding maximal ideal
$\m_h \subset R(Z)$.

The next lemma shows that the action of $R(Z)$ on $G(G,S_\Psi)$
is independent of the choice of $h \in \Psi$. This fact is
crucial for the proof of the nonabelian localization theorem.

\begin{lemma} \label{lem.hindep}
Let $\Psi$ be a semisimple conjugacy class and let $h_1$ and $h_2= kh_1k^{-1}$
be
in $\Psi$. Set $Z_1 = {\mathcal Z}_G(h_1)$ and $Z_2 = {\mathcal Z}_G(h_2)$.  
Let $C_k: Z_1 \rightarrow Z_2$ denote conjugation by $k$, and let
$C_k^*: R(Z_2) \rightarrow R(Z_1)$ denote the pullback.  Then
$C_k^*(\m_{h_2}) = \m_{h_1}$ and the pullback $C_k^*$  is compatible 
with the actions of  $R(Z_1)$ and $R(Z_2)$ on $G(G,S_\Psi)$
as defined above.
\end{lemma}
\begin{proof}
Since $C_k(h_1) = h_2$, we have $C_k^*(\m_{h_2}) = \m_{h_1}$.  To verify the compatibility
of the $R(Z_i)$-actions, we must unwind the definitions.
The maps $a_i: G/Z_i \rightarrow \Psi$ defined by $a_i(g Z_i) = g h_i g^{-1}$ give 
$G$-equivariant identifications of $G/Z_i$ with $\Psi$.  These in turn give change of group identifications
of $R(Z_i)$ with $K_0(G, \Psi)$, and the action of $R(Z_i)$ on $G(G, S_{\Psi})$ is obtained
composing these identifications with the pullback $K_0 (G, \Psi) \rightarrow K_0(G, S_{\Psi})$.  Therefore, it suffices to prove that $C_k^*$ is compatible with the identifications of
$R(Z_i)$ with $K_0(G, \Psi)$.  

The identifications of $G/Z_i$ with $\Psi$ are compatible with the isomorphism
$\kappa: G /Z_1 \rightarrow G/Z_2$ taking $g Z_1$ to $g k^{-1} Z_2$.  If $V_2$ is a $Z_2$-module, write $C_k^* V_2$ for the same underlying vector space, but with $Z_1$-module structure obtained by pullback via the map $C_k$.  The identification of $R(Z_i)$ with $K_0(G, \Psi)$ takes the class of a $Z_i$-module $V_i$ to the class of the vector bundle $G \times_{Z_i} V_i$ on $G/Z_i$.  The pullback via $\kappa$ of the vector bundle
$G \times_{Z_2} V_2$ is $G$-equivariantly isomorphic to the vector bundle $G \times_{Z_1} C_k^* V_2$, which proves the result.
\end{proof}

\begin{prop} \label{prop.nfinvertible}
Let $G$ be an algebraic group acting 
on a smooth algebraic space $X$. Let $\Psi = C_G(h)$ be 
a semi-simple conjugacy class, let $Z = {\mathcal Z}_G(h)$, and let $f \colon S_\Psi \to X$
be the projection. Let 
$N_f^* = [T_{S_{\Psi}}^*] - [f^*(T_X^*)]  \in K_0(G,S_\Psi)$
be the class of the relative cotangent bundle.
The map $\cap\; \lambda_{-1}(N_f^*)\colon G(G,S_\Psi) \to
G(G,S_\Psi)$ restricts to an isomorphism after localizing
at the maximal ideal $\m_h \in R(Z)$.
\end{prop}
\begin{remark} \label{rem.nfinvertible}
Recall that $S_\Psi$ is smooth when $X$ is smooth (Remark \ref{rem.spsismooth})
so that $f \colon S_\Psi \to X$ is a morphism of smooth algebraic spaces.
\end{remark}
\begin{proof}
Let $\iota\colon X^h \to X$ be the 
inclusion of the fixed locus
of $h$. By Theorem \ref{thm.weakloc} the map
$$\cap \; \lambda_{-1}(N^*_\iota) \colon G(Z,X^h)_{\m_{h}} 
\to G(Z,X^h)_{\m_{h}}$$
is an isomorphism. Let $i \colon S_\Psi \to G \times_Z X$
be the corresponding inclusion of mixed spaces.
By Morita equivalence 
the map
$$\cap \; \lambda_{-1}(N_i^*) \colon G(G,S_\Psi)_{\m_h} \to
G(G,S_\Psi)_{\m_h}$$ is an isomorphism.

Let $\eta \colon G \times_Z X \to X$ be the projection. Then
$f = \eta \circ i$, so $[T_f] = [N_i] - [i^*T_\eta] \in K_0(G,S_\Psi)$
Thus $\lambda_{-1}(N_i^*) = \lambda_{-1}(N_f^*) \lambda_{-1}(i^*T_\eta^*)$.
Since $\lambda_{-1}(N_i^*)$ acts by automorphisms on $G(G,S_\Psi)_{\m_h}$ 
it follows
that $\lambda_{-1}(N_f^*)$ must as well.
\end{proof}

Now assume that $G$ is connected and that $Z = {\mathcal Z}_G(h)$ is 
a Levi
factor of a parabolic subgroup $P \subset G$. Since $P/Z$ is isomorphic
to affine space we may identify identify $G(G, G\times_P X)$ with
$G(Z,X)$.

\begin{prop} \label{prop.plocal}
Let $G$ be a connected algebraic group acting on a smooth algebraic space $X$.   Let $\Psi = C_G(h)$; assume that $Z = {\mathcal Z}_G(h)$ is a Levi factor of a parabolic subgroup $P$.  Let 
$i$ and $j$ be the maps of $S_{\Psi}$ into $G \times_Z X$ and $G \times_P X$ defined above.

(a) The map 
$$\cap \; \lambda_{-1}(N_{j}^*) \colon G(G,S_\Psi) \to G(G,S_\Psi)$$ is an
isomorphism
after localizing at the maximal ideal $\m_h \subset R(Z)$.

(b) If in addition $G$ acts on $X$ with finite stabilizer and if $\beta \in G(G,G \times_P X)_{\m_h}$, then 
\begin{equation} \label{eqn.parabolic} \beta = 
j_* \left(\lambda_{-1}(N_j^*)^{-1} \cap j^* \beta\right).
\end{equation}
where $\lambda_{-1}(N_j^*)^{-1} \cap j^*\beta$ denotes
the image of $j^*\beta$ under the inverse of the map 
$\cap\; \lambda_{-1}(N_j^*)$.
\end{prop}
\begin{proof} 
(a) The map $j$ factors as $\rho \circ i$ where $\rho \colon G \times_P X
to G \times_Z X$ is the projection. Thus
$[N_j] = [N_i] - [i^*T_\rho]$. Hence
$\lambda_{-1}(N_i^*) = \lambda_{-1}(N_j^*)\lambda_{-1}(T_\rho^*)$.
By Morita equivalence and Theorem \ref{thm.weakloc}, 
$\lambda_{-1}(N_i^*)$ acts by automorphisms on
$G(G,S_\Psi)_{\m_h}$. Therefore, $\lambda_{-1}(N_j^*)$ must as well.

(b) By Morita equivalence
and Theorem \ref{thm.weakloc} (for the group $Z$)
the pushforward  $i_*\colon G(G,S_\Psi)_{c_\Psi} \to G(G,G \times_Z X)_{\m_h}$ 
is an isomorphism. Since the action of $\lambda_{-1}(N_i^*)$ on
$G(G,S_\Psi)_{c_\Psi}$
is invertible, it follows from 
the self intersection formula 
that 
$i^*\colon G(G,G \times_Z X)_h \to G(G,S_\Psi)_{c_\Psi}$ is also
an isomorphism. Hence 
$$j^* = (\rho \circ i)^* \colon G(G,G \times_P X)_{\m_h}
\to G(G,S_\Psi)_{c_\Psi}$$ is as well.
Thus, it suffices to prove that Equation \eqref{eqn.parabolic}
holds after applying $j^*$ to both sides.
By the self-intersection formula for the regular
embedding $j$, we have
$$
j^*j_*\left(\lambda_{-1}(N_j^*)^{-1} \cap  j^*\beta \right)
= \lambda_{-1}(N_j^*) \cap \left( \lambda_{-1}(N_j^*)^{-1} \cap j^*\beta\right)
= j^*\beta
$$
\end{proof} 

\subsection{The central summand}
Suppose that $G$ acts with finite stabilizers on $X$.  Let
$\Psi = C_G(h)$, and keep the notation above.
Since $h$ is 
a one-element conjugacy class contained in $\Psi \cap Z$, Proposition 
\ref{prop.decomp} implies that
$G(G,S_\Psi)_{\m_h}$ is a summand
in $G(G,S_\Psi)_{\m_\Psi}$. By Lemma \ref{lem.hindep}, this summand
is independent of choice of $h \in \Psi$.

\begin{definition} \label{def.centralpiece} 
Let $G$ act with finite stabilizers on $X$, and let $\Psi = C_G(h)$.
With notation as above,
the summand  $G(G,S_{\Psi})_{\m_h} \subset G(G,S_{\Psi})_{\m_\Psi}$
(which is independent of $h \in \Psi$) will be called the central summand
and denoted 
$G(G,S_\Psi)_{c_\Psi}$.   
The component of $\beta \in G(G,S_\Psi)$ in this summand will be 
denoted $\beta_{c_\Psi}$.
\end{definition}

\section{The nonabelian localization theorem} \label{s.nonabelian}
The main theorem of our paper is the following nonabelian localization theorem.

\begin{thm}[Explicit nonabelian localization]\label{thm.localization}
Let $G$ be an algebraic group acting with finite stabilizer on
a smooth algebraic space $X$. Let $\Psi = C_G(h)$ be a semisimple
conjugacy class and let $f \colon S_\Psi \to X$ be the projection.
If $\alpha \in G(G,X)$ let $\alpha_\Psi$ be the component
of $\alpha$ supported at the maximal ideal $\m_\Psi \subset R(G)$. Then
\begin{equation} \label{eqn.werock} 
\alpha_\Psi = 
f_* \left( \lambda_{-1}(N_f^*)^{-1} \cap (f^* \alpha)_{c_\Psi}\right)
\end{equation}
where $\lambda_{-1}(N_f^*)^{-1} \cap (f^*\alpha)_{c_\Psi}$
is the image of $(f^*\alpha)_{c_\Psi}$ under the inverse
of the automorphism $\cap \; \lambda_{-1}(N_f^*)$ of $G(G,S_\Psi)_{c_\Psi}$.
\end{thm}

The theorem can be restated in way that is sometimes more
useful for calculations.  As usual, let $Z = {\mathcal Z}_G(h)$.  
Let 
$\iota^! \colon G(G,X) \to G(Z,X^h)$ be the composition of the restriction
functor $G(G,X) \to G(Z,X)$ with the pullback $G(Z,X) \stackrel{\iota^*} \to 
G(Z,X^h)$.  Let $\beta_h$ denote the component of
$\beta \in G(Z,X^h)$ in the summand $G(Z,X^h)_{\m_h}$.
Let ${\mathfrak g}$ (resp.~${\mathfrak z}$) be the
adjoint representation of $G$ (resp.~$Z$). The restriction of
the adjoint representation to the subgroup $Z$ makes ${\mathfrak g}$
a $Z$-module, so ${\mathfrak g}/{\mathfrak z}$ is a $Z$-module.
Let $\eta \colon G \times_Z X \to X$ be the projection.
By Remark \ref{rem.tantomixedbund}, $T_\eta = G \times_Z (X \times
{\mathfrak g}/{\mathfrak z})$. We therefore obtain the following 
corollary.
\begin{cor} \label{cor.localization}
With assumptions as in Theorem \ref{thm.localization},
let $h$ be an element of $\Psi$.  
Let  $Z = {\mathcal Z}_G(h)$, and let
$\iota_! \colon G(Z,X^h) \to G(G,X)$ be the map obtained
by composing $f_*$ with the Morita equivalence isomorphism $G(Z,X^h) \to
G(G, S_\Psi)$.  If
$\alpha \in G(G,X)_{\m_\Psi}$, then
\begin{equation} 
\alpha = \iota_! \left( \lambda_{-1}(N_\iota^*)^{-1} \cap
\lambda_{-1}(({\mathfrak g}/{\mathfrak z})^*)\cap
(\iota^!\alpha)_h\right).
\end{equation}
\end{cor}
\subsection{Proof of Theorem \ref{thm.localization} if $G$ is connected and $Z$ is a Levi factor of a parabolic subgroup}
\label{sect.localizationproof1}
In this section we prove Theorem \ref{thm.localization} under the
assumptions that $G$ is connected and that $Z=
{\mathcal Z}_G(h)$ is a Levi factor
in a parabolic subgroup $P \subset G$.  (In fact, it would suffice for our purposes to
take $G$ equal to a product of general linear groups, and then $Z$ is automatically
such a Levi factor.)   Choose a maximal torus $T$ and
and Borel subgroup $B$ such that $h \in T \subset B \subset P$.
We have maps
of mixed spaces
$$
G \times_B X\stackrel{p}{\rightarrow} G \times_P X\stackrel{q}{\rightarrow} X;
$$
and we write $\pi = q \circ p$.  Let $j: S_{\Psi} \rightarrow G
\times_P X$ denote the regular embedding of Lemma
\ref{lemma.zpsitoxp}.  

Since $G \times_T X \rightarrow G \times_B X$
is a bundle with fibers isomorphic to $B/T$, pullback gives an isomorphism of
$G(G, G \times_B X)$ with $G(G, G \times^T X)$.  Now, $\Psi \cap T$ consists of a finite number
of elements $h= h_1, \ldots , h_w$, where 
$w = |W(G,T)|/|W(Z,T)|$.
By Proposition
\ref{prop.changedecomp}, if $\alpha \in G(G,X)_{\m_{\Psi}}$, then there is a decomposition
$$\pi^*\alpha = \sum_{l = 1}^w (\pi^*\alpha)_{h_l},$$
where $(\pi^*\alpha)_{h_l}$ refers to the component supported
at $\m_{h_l} \subset R(T)$. 
Because $h$ is central in $Z$, there is a also component
$(q^* \alpha)_h$ of $q^* \alpha$ supported at the maximal
ideal $\m_h \subset R(P)= R(Z)$.

The key step is given by the following proposition.

\begin{prop} \label{prop.decomposition}
Keep the assumptions of Theorem \ref{thm.localization}, and in
addition assume that $G$ is connected and $Z$ is a Levi factor of a parabolic subgroup $P$.
Then the following identity holds in $G(G,X)_{\m_\Psi}$:
\begin{equation} 
f_* \left(\lambda_{-1}(N_f^*)^{-1} \cap (f^*\alpha)_{c_{\Psi}} \right) = 
\frac{1}{|W(Z,T)|} \pi_*
\left(\lambda_{-1}(T_\pi^*) \cap (\pi^*\alpha)_h\right).
\end{equation}
\end{prop}
\begin{proof}
We keep the notation introduced before the statement of the
proposition.  Since $f = q \circ j$ and $j^*$ is an $R(Z)$-module homomorphism
$j^*((q^*\alpha)_h) =
(f^*\alpha)_{h}$.  But $(f^*\alpha)_h$ is independent
of $h$ and equals $(f^*\alpha)_{c_\Psi}$
by definition of the central summand. Now, $[N_f^*] = [N_j^*] - [j^*T_q^*]$, so
\begin{equation} \label{eqn.firstlambda}
\begin{array}{rcl}
\lambda_{-1}(N_f^*)^{-1} \cap (f^*\alpha)_{c_{\Psi}} &  = & \lambda_{-1}(N_j^*)^{-1}
\cap (\lambda_{-1}(j^*T_q^*) \cap (f^*\alpha)_{c_{\Psi}})\\
& = & 
\lambda_{-1}(N_j^*)^{-1} \cap j^*\left(\lambda_{-1}(T_q^*) \cap (q^*\alpha)_h
\right).
\end{array}
\end{equation}
Applying $f_* = q_* \circ j_*$ to both sides of Equation
\eqref{eqn.firstlambda}, we obtain
\begin{equation} \label{eqn.secondlambda}
\begin{array}{rcl}
f_*\left(\lambda_{-1}(N_f^*)^{-1} \cap (f^*\alpha)_{c_\Psi}\right)
& = & q_* j_* \left(\lambda_{1}(N_j^*)^{-1} \cap j^*(\lambda_{-1}(T_q^*)
\cap (q^*\alpha)_h)\right)\\
& = & q_*\left(\lambda_{-1}(T_q^*) \cap (q^*\alpha)_h\right)
\end{array}
\end{equation}
where the second equality follows from Proposition \ref{prop.plocal}.
By compatibility of pullback with support, $(\pi^*\alpha)_h = p^*((q^*\alpha)_h)$.
By Proposition \ref{prop.tanpush},
$$\lambda_{-1}(T_q^*) \cap (q^*\alpha)_h = \frac{1}{|W(Z,T)|}
p_*\left(\lambda_{-1}(T_\pi^*) \cap (\pi^*\alpha)_h\right)$$ 
so 
$$q_*\left(\lambda_{-1}(T_q^*) \cap (q^*\alpha)_h\right)
= \frac{1}{|W(Z,T)|}\pi_*\left(\lambda_{-1}(T_\pi^*) \cap 
(\pi^*\alpha)_h\right).$$
\end{proof}

The proof of the theorem for $G$ connected and $Z$ a Levi factor of a parabolic subgroup is an easy consequence
of the preceding proposition.
By Proposition \ref{prop.tanpush},  
$$\alpha = \frac{1}{|W(G,T)|} \pi_*\left(
\lambda_{-1}(T_\pi^*) \cap  \pi^* \alpha\right).$$
Therefore, 
\begin{eqnarray*}
\alpha & = & \frac{1}{|W(G,T)|}\sum_{l = 1}^w \pi_*\left(
\lambda_{-1}(T_\pi^*) \cap  (\pi^* \alpha)_{h_l}\right)\\
& = &  
 \frac{|W(Z,T)|}{|W(G,T)|}\sum_{l = 1}^w
f_*\left( \lambda_{-1}(N_f^*)^{-1} \cap (f^*\alpha)_{c_\Psi}\right)\\
& = & f_*\left( \lambda_{-1}(N_f^*) \cap (f^*\alpha)_{c_\Psi}\right),
\end{eqnarray*}
as desired.

\begin{remark}
As a corollary of the proof in this case we obtain
the following induction formula for $\alpha \in G(G,X)_{\m_\Psi}$. 
\begin{equation}
\alpha = \frac{|W(G,T)|}{|W(Z,T)|}\pi_*\left(\lambda_{-1}(T_\pi^*) \cap (\pi^*\alpha)_h\right)
\end{equation}
\end{remark}

\subsection{Proof of Theorem \ref{thm.localization} for arbitrary $G$}\label{sect.localizationproof2}
By Proposition \ref{prop.goodembedding}, we can embed $G$ as a closed
subgroup of a product of general linear groups $Q$ such that, writing
$\Psi = C_G(h)$ and
$\Psi_Q = C_Q(h)$, we have $\Psi_Q \cap G = \Psi$.  
Write $Z = {\mathcal Z}_G(h)$ and $Z_Q = {\mathcal Z}_Q(h)$.
Let $Y = Q
\times_G X$; then $Q$ acts with finite stabilizer on $Y$.
The group $Q$ is connected, and direct calculation shows that
$Z_Q$ is a
a product of general linear groups which is a Levi factor of a parabolic. Therefore, 
the nonabelian localization theorem applies
to the $Q$-action on $Y$.

By Proposition \ref{prop.changedecomp} and the fact that $\Psi_Q \cap
G = \Psi$, we have $G(Q,Y)_{\m_{\Psi_Q}} = G(Q,Y)_{\m_\Psi}$.  As usual, $S_{\Psi_Q}$
denotes the part of the global stabilizer $S_Y$ corresponding to $\Psi_Q$.
For
simplicity of notation we will write $S_Q$ for $S_{\Psi_Q}$ and  $S$ for
$S_{\Psi}$.  We denote by $f_Q$ the morphism $S_Q \rightarrow
Y$. If $\beta \in G(Q,S_Q)$ then by definition,
$\beta_{c_{\Psi_Q}}$ is the component of $\beta$ supported at the maximal ideal
$\m_h^Q \subset R(Z_Q)$. Likewise, if $\gamma \in G(Q, Q \times_G S)=
G(G,S)$ then $\gamma_{c_\Psi}$ is the component supported at $\m_h \subset
R(Z)$.  Since $h$ is central in $Z$ and $Z_Q$,  Remark \ref{rem.pickbill}
implies that if $M$ is any
$R(Z)$-module then $M_{\m_h} = M_{\m_h^Q}$, where the $R(Z_Q)$ action
on $M$ is given by the restriction homomorphism $R(Z_Q) \to R(Z)$.
Thus $\gamma_{c_\Psi}$ may also be identified with the component of 
$\gamma$ supported at $\m_h^Q$.

Note that $R(Z_Q)$ acts on $G(G,S_Q)$ and $R(Z)$ acts on
$G(Q, Q \times_G S) = G(G,S)$.  By Morita equivalence,
the localization theorem for $G$
acting on $X$ is a consequence of the following lemma.

\begin{lemma}
Keep the assumptions and notation of this subsection.  Then there is a
$Q$-equivariant isomorphism $\Phi: Q \times_G S \to S_Q$
such that:

i) $\Phi^* \colon G(Q,S_Q) \to G(Q, Q \times_G S)$ is an $R(Z_Q)$-module
homomorphism where the action of $R(Z_Q)$ on $G(Q,Q \times_G S)$ is given
by the restriction homomorphism $R(Z_Q) \to R(Z)$.

ii) $f_Q \circ \Phi = (1 \times_G f)$.
\end{lemma}

\begin{proof}
Consider the map 
$$T \colon Q  \times S \to \Psi_Q \times Q \times X, \; \; (q,g,x) \mapsto
(qgq^{-1},q,x).$$
This map induces a map of quotient spaces 
$\tilde{\Phi} \colon Q \times_G S \to \Psi_Q  \times Y$ such that the following diagram commutes:
\begin{equation} \label{diag.qhom}
\begin{array}{ccc} Q \times_G S & \to & \Psi_Q \times Y\\
\downarrow & & \downarrow \\
Q \times_G \Psi & \stackrel{\phi}{\to} & \Psi_Q
\end{array}
\end{equation}
Here the vertical arrows are the obvious projections and
the bottom horizontal arrow is given by $\phi(\overline{(q,k)} ) =  qkq^{-1}$.
Note that $S_Q$ is a closed subspace of $\Psi_Q \times Y$.
We have a commutative diagram
$$
\begin{array}{ccccc}
R(Z_Q) & \stackrel{=}\to & K_0(Q,\Psi_Q) & \to & K_0(Q,\Psi_Q \times Y) \\
\downarrow & & \downarrow \phi^* &  & \downarrow \tilde{\Phi}^* \\
R(Z) & \stackrel{=}\to & K_0(Q,  Q \times_G  \Psi  ) & \to & K_0(Q,Q \times_G S).
\end{array}
$$
Here the second arrow in each row is a pullback map.  The first arrow
in the top row takes $[V] \in R(Z_Q)$ to the class of the vector
bundle $Q \times_{Z_Q} V$ on $Q/Z_Q = \Psi_Q$.  The first arrow in the
bottom row takes $[W] \in R(Z)$ to the class of the vector bundle $Q
\times_Z W$ on $Q/Z = Q \times_G \Psi$.  The action of $R(Z_Q)$ on
$G(Q,S_Q)$ (resp.~of $R(Z)$ on $G(Q,Q \times_G S)$) is defined using
the composition along the top (resp.~bottom) row.  The commutativity
of this diagram implies that $\tilde{\Phi}^*$ is a $R(Z_Q)$-module
homomorphism.

Next we show that $\tilde{\Phi}$ induces an isomorphism of $Q \times_G S$ onto
the subspace $S_Q \subset \Psi_Q \times Y$.  Let $W$
be the inverse image of $S_Q$ in $\Psi_Q \times Q \times X$.
Then $W$ is the closed subspace consisting 
of triples $(k,q, x)$ such that $g = q^{-1}kq$ is in
$G$ and $g x = x$.  If $(g,x) \in S$ and $q \in Q$ then clearly
$(qgq^{-1},q,x) \in W$.  Thus $\tilde{\Phi}$ factors through a
morphism $\Phi \colon Q \times_G S \to S_Q$.  Since $\tilde{\Phi}^*$ is
an $R(Z_Q)$-module homomorphism, so is $\Phi^*$.  
Next we show that $\Phi$ is an isomorphism. 
Since we work over $\C$
and $\Phi$ is a representable morphism of smooth (hence normal)
algebraic spaces, it suffices, by Zariski's main theorem (proved for
algebraic spaces in \cite[Theorem V4.2]{Knu:71}), to prove that $\Phi$
is bijective on geometric points.

First we show that $\Phi$ is injective.
Let $(q_1,g_1,x_1)$ and $(q_2,g_2,x_2)$ be two points
of $Q \times S$ such that $T(q_1,g_1,x_1)= (q_1g_1q_1^{-1},q_1,x_1)$ and 
$T(q_2,g_2,x_2) = (q_2g_2q_2^{-1},q_2,x_2) $ have
the same image in $\Psi_Q \times Y$. Then $q_2g_2q_2^{-1} = q_1g_1q_1^{-1}$
and there is an element
$g \in G$ such that $x_2 = gx_1$ and $q_2 =  q_1 g^{-1}$. Thus
$g_2 = gg_1g^{-1}$ and hence $(q_1,g_1,x_1)$ and $(q_2, g_2, x_2)$
are in the same $G$-orbit in $Q \times S$. Therefore, $\Phi$ is injective
on geometric points.

Conversely, suppose that $(k,q,x) \in W$. Then 
$g =q^{-1}kq \in \Psi_Q \cap G = \Psi$. Thus $(k,q,x) = T(q,g,x)$
so $\Phi$ is surjective on geometric points.

Finally, the fact that $f_Q \circ \Phi = (1 \times_G f)$
is clear from the definition of $\Phi$.
This completes the proof of the lemma and with it
Theorem \ref{thm.localization}.
\end{proof} 
\section{Riemann-Roch for quotients} \label{s.err}
As an application of the nonabelian localization theorem, we can give
an explicit formula for the Riemann-Roch map for quotients of smooth
algebraic spaces by proper actions of algebraic groups. Recall,
(Section \ref{ss.group}) that any algebraic space with a proper
group action is automatically separated.

Before
stating this, we need some preliminary results about invariants.

\subsection{Invariants and equivariant $K$-theory}
Let $G$ be an algebraic group acting properly on an algebraic
space $X$. The theorem of Keel and Mori \cite{KeMo:97} states that the
quotient stack $[X/G]$ has a moduli space $Y = X/G$ in the category of
algebraic spaces. Translated in terms of group actions, this means that
the map of algebraic spaces $X \stackrel{\pi} \to Y$ is a categorical
geometric quotient in the category of algebraic spaces.

\begin{definition}
Let $X$ be an algebraic space with a proper $G$ action
and let $X \stackrel{\pi} \to Y$ be the geometric quotient.
If ${\mathcal F}$ is a $G$-equivariant quasi-coherent ${\mathcal O}_X$-module,
define ${\mathcal F}^G  = (\pi_*{\mathcal F})^G$.  We call this the
functor of taking $G$-invariants.
\end{definition}

To obtain a map from equivariant $K$-theory to the $K$-theory of the quotient,
we need the following lemma.

\begin{lemma} \label{lem.invariants} 
Let $X$ be an algebraic space with a proper $G$-action (see Section
\ref{ss.group}) and let $X \stackrel{\pi} \to Y$ be the geometric
quotient.  The assignment ${\mathcal F} \mapsto {\mathcal F}^G$
defines an exact functor ${\tt coh}^G_X \to {\tt coh}_Y$.
\end{lemma}

\begin{remark} If $X$ and $Y$ are both schemes, and $G$ is reductive,
then this lemma is a consequence of the facts that the quotient map $X
\rightarrow Y$ is affine \cite[Proposition 0.7]{MFK:94}, and that
taking invariants by a locally finite action of a reductive
group is an exact functor (in characteristic $0$).
\end{remark}

The theorem of Keel and Mori is proved 
using an \'etale local
description of the moduli space $Y = X/G$. In particular they prove
that $[X/G]$ has a representable \'etale cover by quotient
stacks $\{ [U_i/H_i]\}$ with $U_i$ affine, $H_i$ finite and such that
the following diagram of stacks and moduli spaces is Cartesian:
$$\begin{array}{ccc}
[U_i/H_i] & \to & [X/G] \\
\downarrow & & \downarrow\\
U_i/H_i & \to & X/G
\end{array}$$
Let $V_i = U_i \times_{[X/G]} X$.
Then $V_i$ is affine and has commuting and free actions of $G$ and $H_i$. 
Let $X_i = Y_i/H_i$. Since $H_i$ acts freely, $X_i \to X$ is \'etale.
If we let $Y_i = X_i/G_i$, then the map of quotients
$Y_i \to Y$ is also \'etale, and the following diagram of spaces and quotients
is Cartesian:
\begin{equation} \label{diag.structure}
\begin{array}{ccc}
X_i & \to & X\\
\downarrow & & \downarrow\\
Y_i & \to & Y
\end{array}
\end{equation}
The actions of $G$ and $H_i$ on $Y_i$ commute so
$X_i/G = (Y_i/G_i)/H_i$.
The local structure of geometric quotients given by Diagram 
\eqref{diag.structure} implies that the quotient map $X \to Y$ is
an affine morphism in the category of algebraic spaces.

\begin{proof}[Proof of Lemma \ref{lem.invariants}]
The question is local in the \'etale topology on the quotient $Y$.
Thus we may assume that there is an affine scheme $V$ and a finite
group $H$ such that $V$ has commuting free actions by $H$ and $G$
and $X = V/G$.
It follows that the $G \times H$-quotient map $V \stackrel{p} \to Y$
factors as $ p = q\circ \pi_V = \pi \circ q_V$ where
$\pi_V \colon V \to V/G$ is a $G$-torsor, $q_V \colon V \to X$
is an $H$-torsor and $q\colon V/G \to Y$ is a quotient by $H$.

Let ${\mathcal G}$ be a $(G \times H)$-equivariant coherent sheaf
on $V$. Then 
\begin{eqnarray*}
{\mathcal G}^{G \times H}  & = & ({\mathcal G}^G)^H \\
& = & ({\mathcal G}^H)^G
\end{eqnarray*}
The group $H$ is finite and $Y = (V/G)/H$ so  
the assignment ${\mathcal H} \mapsto {\mathcal H}^H$
is an exact functor ${\tt coh}_{V/G}^H \to {\tt coh}_Y$.
Since $G$ and $H$ act freely on $V$ the assignments
${\mathcal G} \mapsto {\mathcal G}^G$ and ${\mathcal G} \mapsto
{\mathcal G}^H$ define equivalences ${\tt coh}^{G \times H}_V \to {\tt
coh}^H_{V/G}$ and ${\tt coh}^{G \times H}_V \to {\tt coh}^G_X$.
Thus, the assignment ${\mathcal G} \mapsto {\mathcal G}^{G \times H}$
is an exact functor ${\tt coh}^{G \times H}_V \to 
{\tt coh}_X$. Since the assignment ${\mathcal G} \mapsto {\mathcal G}^{H}$
is an equivalence ${\tt coh}^{G \times H}_V \to {\tt coh}^G_X$
it follows that the assignment ${\mathcal F} \to {\mathcal F}^G$
is an exact functor ${\tt coh}^G_X \to {\tt coh}_Y$.
\end{proof}

If $X$ is an algebraic space with a proper $G$-action 
and geometric quotient $X \stackrel{\pi} \to Y$ let
$\pi_G \colon G(G,X) \to G(Y)$ be the map on $K$-theory
induced by the exact functor 
${\tt coh}_X^G \to {\tt coh}_Y$ given by 
${\mathcal F} \mapsto {\mathcal F}^G$.  We will usually denote
$\pi_G(\alpha)$ by $\alpha^G$.

If $G$ acts properly on $X$ and $q \colon X' \to X$ is a finite
$G$-equivariant map then $G$ acts properly on $X'$ \cite[Prop. 2.1]{EdGr:03}.
Let $\pi \colon X \to Y$, $\pi'\colon X' \to Y'$ be the geometric quotients.
Since the composite map $\pi \circ q \colon X' \to Y$ is $G$-invariant
there is a map of quotient $q'\colon Y' \to Y$ such that the diagram commutes.
\begin{equation} \label{diag.invcommute}
\begin{array}{ccc}
X' & \stackrel{q} \to & X \\
\pi^{'} \downarrow & & \pi \downarrow\\
Y^{'} & \stackrel{q^{'}} \to & Y
\end{array}
\end{equation}

\begin{lemma} \label{lem.invcommute}
The map $q'$ is finite and $\pi_G \circ q_*  = q'_*  \circ \pi'_G $
as maps $G(G,X') \to G(Y)$.
\end{lemma}
\begin{proof}
Working locally in the \'etale topology we may assume that $Y$ is
affine. It follows (since $q$ is finite) that all of the other spaces
in Diagram \eqref{diag.invcommute} are affine.  Let $X = \Spec A$, $X'
= \Spec B$. Then $Y' = \Spec B^G$ and $Y = \Spec A^G$. If $M$ is a
$G$-equivariant $B$-module then $_{A^G} M^G = (_A M)^G$
(c.f. \cite[Proposition 2.3]{EdGr:03}).  Translated to sheaves this
means that ${\mathcal F}$ is a $G$-equivariant quasi-coherent sheaf on
$X'$ then $(q_* {\mathcal F})^G = q'_*({\mathcal F}^G)$ as
quasi-coherent sheaves on $Y$.

To complete the proof of the Lemma we need to show that $q'$ is
finite. Since $q'$ is affine we need to show that $q'_* {\mathcal
O}_{Y'}$ is coherent. Now $\pi'$ is a geometric quotient so ${\mathcal
O}_{Y'} = {\mathcal O_{X'}^G}$.  Thus $q_*{\mathcal O}_{Y'} = (q_*
{\mathcal O}_{X'})^G$ is coherent because $q$ is finite.
\end{proof}

\subsection{Twisting equivariant $K$-theory by a central subgroup}

Let $Z$ be an algebraic group and let $H$ be a subgroup (not
necessarily closed) of the center of $Z$ consisting of semisimple
elements. If $V$ is any representation of $Z$ then, since any
commuting family of semisimple endomorphisms is simultaneously
diagonalizable, $V$ can be written as a direct sum of $H$-eigenspaces $V
= \oplus V_{\chi}$ where the sum is over all characters $\chi: H
\rightarrow \C^*$ and $H$ acts on $V_{\chi}$ by $h\cdot v = \chi(h)v$.
Since elements of $H$ are central in $Z$, each eigenspace $V_{\chi}$
is $Z$-stable, and therefore we have $[V] = \sum [V_{\chi}]$ in
$R(Z)$.  We define an action of $H$ on the representation ring $R(Z)$,
by
$$
h \cdot [V] = \sum \chi (h)^{-1} [V_{\chi}],
$$
for $h \in H$ and $[V] = \sum [V_{\chi}] \in R(Z)$.  Because $H$ is a
central subgroup, it acts on $\C[Z]^Z$ by the rule $(h \cdot f)(z)
=f(h^{-1}z)$, for $h \in H$, $z \in Z$, $f \in \C[Z]^Z$.  With these
actions, the character map from $R(Z)$ to $\C[Z]^Z$ is $H$-equivariant.
Recall that maximal ideals of $R(Z)$ are of the form $\m_{\Psi}$,
where $\Psi$ is a semisimple conjugacy class in $Z$.  The
$H$-equivariance of the character map implies that if $h \in H$, then
$$
h \cdot \m_{\Psi} = \m_{h \Psi}.
$$
In particular, 
$$
h^{-1} \m_h = \m_1,
$$
the augmentation ideal of $R(Z)$. 

More generally, let $X$ be an algebraic space with a $Z$-action, such
that $H$ acts trivially on $X$.  In this case, if $\E$ is any
$Z$-equivariant coherent sheaf on $X$, and $\chi$ a character of $H$,
let $\E_{\chi}$ be the subsheaf of $\E$ whose sections (on any \'etale
open set $U$) are given by
$$
\E_{\chi} (U) = \{ s \in \E(U) \ | \ h \cdot s = \chi (h) s \}.
$$
Then $\E = \oplus_{\chi} \E_{\chi}$ is a decomposition of $\E$ into a
direct sum of $H$-eigensheaves for $\E$.  We define an action of $H$ on
$G_0(Z,X)$ by
$$
h \cdot [\E] = \sum \chi(h)^{-1} [\E_{\chi}].
$$
If $X$ is a point, this reduces to the previous definition of the
action of $H$ on $R(Z)$.  Thus, the actions of $H$ on $R(Z)$ and on
$G_0(Z,X)$ are compatible: if $r \in R(Z)$, $\alpha \in G_0(Z,X)$, then $h
\cdot (r \alpha) = (h \cdot r) (h \cdot \alpha)$.  If $Z$ acts on $X$
with finite stabilizers, then we can identify $G_0(Z,X) = \oplus
G_0(Z,X)_{\m_{\Psi}}$, where the sum is over finitely many conjugacy
classes $\Psi$.  The compatibility of the actions implies that
\begin{equation} \label{eqn.shift}
h^{-1} G_0(Z,X)_{\m_h} = G_0(Z,X)_{\m_1}.
\end{equation}
We refer to the action of $H$ on $G_0(Z,X)$ as twisting.  If $\alpha \in G_0(Z,X)$, we will often write
$\alpha(h)$ for $h^{-1} \cdot \alpha$ and refer to this map as twisting by $h$.

\subsection{Twisting and the global stabilizer}
We now define a central twist.  Let $\Psi = C_G(h)$ be a semisimple
conjugacy class and let $Z = {\mathcal Z}_G(h)$ and let $H$ be the
cyclic subgroup of $Z$ generated by $h$. Since $h$ is central in $Z$,
there is an action of $H$ on $G_0(Z,X^h)$. The map $\Phi_h \colon G
\times X^h \to S_\Psi$, $(g,x)\mapsto (ghg^{-1},gx)$ identifies
$S_\Psi$ with $G \times_Z X$. Thus by Morita equivalence there is an
action of $H$ on $G_0(G,S_\Psi)$.

We define a twisting map
$$
\begin{array}{c}
G_0(G,S_{\Psi}) \rightarrow G_0(G,S_{\Psi}) \\
\alpha \rightarrow \alpha(c_{\Psi})
\end{array}
$$
to be the map induced by Morita equivalence corresponding to the 
endomorphism of $G_0(Z,X^h)$ given by twisting by $h$.  
We will call this map the central twist.

\begin{lemma}  \label{lemma.independenttwist}
With notation as above, the map $\alpha \mapsto \alpha(c_{\Psi})$ is
independent of the choice of $h \in \Psi$.  Moreover, if
$\alpha \in G_0(G, S_{\Psi})_{c_{\Psi}}$, then $\alpha(c_{\Psi}) 
\in G_0(G,S_{\Psi})_{\m_1}$.  For arbitrary $\beta \in G_0(G,S_{\Psi})$,
the component of $\beta(c_{\Psi})$ in $G_0(G,S_{\Psi})_{\m_1}$
equals $\beta_{c_{\Psi}}(c_{\Psi})$.
\end{lemma}

\begin{proof}
Suppose that $h_1$ and $h_2 = k h_1 k^{-1}$ are two elements of
$\Psi$; let $Z_i = {\mathcal Z}_G(h_i)$.  Then $k X^{h_1} = X^{h_2}$ and $k Z_1
k^{-1} = Z_2$.  We have inclusions $f_i: X^{h_i} \rightarrow S_{\Psi}$
given by $f_i(x) = (h_i, x)$ for $x \in X^{h_i}$.  The equivalence of
categories between ${\tt coh}^G_{S_{\Psi}}$ and
${\tt coh}^{Z_i}_{X^{h_i}}$ takes a
 $G$-equivariant coherent
sheaf $\F$ on $S_{\Psi}$ to the sheaf $f_i^* \F$ on $X^{h_i}$ (since
$f_i$ is $Z_i$-equivariant, the pullback $f_i^* \F$ has a natural
structure of $Z_i$-equivariant sheaf).  Now, $k^*$ induces an
equivalence of categories between ${\tt coh}^{Z_2}_{X^{h_2}}$ and
${\tt coh}^{Z_1}_{X^{h_1}}$, 
and
moreover $k^* f_2^* = f_1^*$.  Therefore, to show the independence of
the twisting map, it suffices to show that a $Z_2$-equivariant sheaf
$\E$ on $X^{h_2}$ is an $h_2$-eigensheaf with eigenvalue $\chi$ if and
only if the pullback sheaf $k^* \E$ on $X^{h_1}$ is an
$h_1$-eigensheaf with eigenvalue $\chi$.  This holds because if $U_2 \to X_2$
is any \'etale open set, and $U_1 = k^* U_1$, then $(k^*\E)(U_1) = \E
(U_2)$, and the automorphism of $(k^*\E)(U_1)$ coming from $h_1$
coincides with the automorphism of $\E(U_2)$ induced by $h_2$.

Next, if $\alpha \in G_0(G, S_{\Psi})_{c_{\Psi}}$, then under the Morita
equivalence isomorphism, $\alpha$ corresponds to an element in 
$G_0(Z,X^h)_{\m_h}$.  By definition of the central twist
and \eqref{eqn.shift}, $\alpha(c_{\Psi})$ corresponds to an element of
$G_0(Z,X^h)_{\m_1}$.  By Proposition \ref{prop.changedecomp},
the Morita equivalence isomorphism identifies
$G_0(Z,X^h)_{\m_1}$ with $G(G_0(G,S_{\Psi})_{\m_1}$.  Hence
$\alpha(c_{\Psi}) \in G_0(G,S_{\Psi})_{\m_1}$.

Finally, if $\beta$ is an arbitrary element of $G_0(G, S_{\Psi})$, then
under the Morita equivalence isomorphism the elements $\beta(c_{\Psi})$ and $\beta_{c_{\Psi}}(c_{\Psi})$
correspond to elements of $G_0(Z,X^h)$ whose components at $\m_1$ are equal,
so the result follows.
\end{proof}

Now assume that $G$ acts properly on $X$. Then it also acts properly
on $S_\Psi$, since the map $f\colon S_\Psi \to X$ is finite.
By Lemma \ref{lem.invariants} 
there is a map in $K$-theory (the map of taking $G$-invariants)
$G_0(G,S_{\Psi}) \rightarrow G_0(S_{\Psi}/G)$, $\alpha \mapsto \alpha^G$.

\begin{lemma} \label{l.ctwistinv}
If $\alpha \in G_0(G,S_{\Psi})$, then $\alpha^G = (\alpha(c_{\Psi}))^G$.
\end{lemma}

\begin{proof}
Let $h \in \Psi$ and $Z = {\mathcal Z}_G(h)$.  Since the map $\Phi_h: G \times
X^h \rightarrow S_{\Psi}$ is a $Z$-torsor, the quotients $X^h/Z$ and
$S_{\Psi}/G$ are both identified with the quotient $M = (G \times X^h)/(G
\times Z)$.  We have a map $G(Z, X^h)
\rightarrow G(M)$ (the map of taking $Z$-invariants) and a map
$G(G,S_{\Psi}) \rightarrow G(M)$ (the map of taking $G$-invariants).
By Morita equivalence, both $G(Z,X^h)$ and $G(G,S_{\Psi})$ are
isomorphic to $G(G \times Z, G \times X^h)$ and under these
isomorphisms, both maps coincide with the map $G(G \times Z, G \times
X^h) \rightarrow G(M)$ (the map of taking $G \times Z$-invariants).

In view of the definition of the central twist, it suffices to prove
that if $\alpha \in G_0(Z,X^h)$ then $\alpha^Z = (\alpha(h))^Z$.  For
this we may assume that $\alpha = [\E]$, where $\E$ is a
$Z$-equivariant coherent sheaf on $X^h$.  As above, write $\E = \oplus
\E_{\chi}$.  Note that $[\E_{\chi}]^Z = 0$ unless $\chi$ is the
trivial character; so, denoting the trivial character by $1$, we have
$\alpha^Z = \sum [\E_{\chi}]^Z = [\E_1]^Z$ and $(\alpha(h))^Z = \sum
\chi(h)^{-1} [\E_{\chi}]^Z = [\E_1]^Z$, completing the proof.
\end{proof}

\subsection{Riemann-Roch for quotients: statement and proof}
If $X$ is a $G$-space let $CH^*_G(X) = \oplus_i CH^i_G(X) \otimes \C$ where
$CH^i_G(X)$ denotes the ``codimension'' $i$ equivariant
Chow groups of $X$ as in \cite[p.~569]{EdGr:00}.
Recall \cite[Theorem 3]{EdGr:98} that if $G$ acts properly on $X$, with quotient $X/G$, then there is an isomorphism
$\phi^G_X: CH^*_G(X) \rightarrow CH^*(X/G)$,
where $CH^*(X/G) = \oplus_i CH^i(X/G) \otimes \C$. The map
is defined as follows: Because $G$ acts with finite stabilizers
$CH^*_G(X)$ is generated by fundamental classes of $G$-invariant cycles.
If $V$ is a closed $G$-invariant subspace let $[V/G]$ be
the image of $V$ under the quotient map. 
Then $\phi^G_X([V]) =\frac{1}{e_V}[V/G]$
where $e_V$ is the order of the stabilizer of a general point of
$V$.

If $G$ and $X$ are understood, we may write
simply $\phi$ or $\phi_X$ for $\phi^G_X$.

If $Y$ is an algebraic space let $\tau_Y \colon G_0(Y) \to CH^*(Y)$ be
the Todd isomorphism of \cite[Theorem 18.3]{Ful:84} (as extended to
algebraic spaces). Likewise if $X$ is a $G$-space let $\tau_X^G \colon
G_0(G,X) \to \Pi_{i=0}^\infty CH^i_G(X)$ be the equivariant Todd map
of \cite{EdGr:00}. When $G$ acts with finite stabilizers, $CH^i_G(X) =
0$ for $i > \dim X$ so the target of the equivariant Todd map is
$CH^*_G(X)$ (cf. \cite[Cor 5.1]{EdGr:00}). Note that $\tau^G_X$ is an
isomorphism only when $G$ acts freely. However, it factors through an
isomorphism $\widehat{G_0(G,X)} \to \Pi_{i = 0}^\infty CH^i_G(X)$
where $\widehat{G_0(G,X)}$ is the completion of $G_0(G,X)$ at the
augmentation ideal of $R(G)$ \cite[Theorem 4.1]{EdGr:00}.

\begin{thm} \label{t.rrtheorem1}
Let $G$ be an algebraic group acting properly on 
a smooth
algebraic space $X$, and let $Y = X/G$ be the quotient.  Fix a
conjugacy class $\Psi$ in $G$ and let $S_{\Psi}$ be the corresponding
part of the global stabilizer, so we have a $G$-equivariant map $f:
S_{\Psi} \rightarrow X$ and an induced map $g: S_{\Psi}/G \rightarrow
Y$ on the quotients.  Let $\alpha \in G_0(G,X)$ and let
$\alpha_{\Psi}$ denote the part of $\alpha$ in $G_0(G,X)_{\m_{\Psi}}$.
Then
\begin{equation} \label{e.rrtheorem1}
\begin{array}{ccl} 
\tau_Y((\alpha_{\Psi})^G) &= &  \phi_G^X\circ \tau^G_X \circ f_*
\left((\lambda_{-1}(N_f^*)^{-1} \cap f^*\alpha)(c_\Psi)\right) \\\\
& = & g_* \circ \phi^G_{S_{\Psi}} \circ \tau^G_{S_{\Psi}}
\left( (\lambda_{-1}(N_f^*)^{-1} \cap f^* \alpha)(c_\Psi) \right).
\end{array}
\end{equation}
\end{thm}

This theorem can be stated in the following equivalent form, which can be more
convenient in applications.

\begin{thm} \label{t.rrtheorem2}
Keep the notation and hypotheses of the previous theorem, and in
addition, fix $h \in \Psi$, let $Z = {\mathcal Z}_G(h)$, and let $\iota:
X^h \rightarrow X$ denote the closed embedding.  Identify $X^h/Z$ with
$S_{\Psi}/G$.  Let $\alpha |_{Z}$ denote the image of $\alpha$ under
the natural map of $G$-equivariant $K$-theory to $Z$-equivariant
$K$-theory.  Then
\begin{equation} \label{e.rrtheorem2}
\tau_Y((\alpha_{\Psi})^G) = g_* \circ \phi^Z_{X^h} \circ \tau^Z_{X^h}
\left( \lambda_{-1}(N^*_{\iota})^{-1} \cap \lambda_{-1}(\g^*/ \z^*)
\cap \iota^*(\alpha|_Z)(h) \right) .
\end{equation}
\end{thm}

\begin{remark}
Note that we do not need to compute $\alpha_{\Psi}$ to apply the
formulas of Theorems \ref{t.rrtheorem1} and \ref{t.rrtheorem2}.
However, the answer would be the same if on the right side of these
formulas we replaced $\alpha$ with $\alpha_{\Psi}$.  The reason is
that the part of the formula corresponding to $\alpha_{\Psi'}$ for any
$\Psi'  \neq \Psi$ vanishes.  Also, since the invariant map is linear,
$\tau_Y(\alpha^G) = \sum_{\Psi} \tau_Y((\alpha_{\Psi}^G)$ can be
computed using the formulas of Theorems \ref{t.rrtheorem1} and
\ref{t.rrtheorem2}.
\end{remark}

\begin{remark}
If $X$ (and hence also $X^h$) is a smooth scheme, then 
$K_0(G,X) = G_0(G,X)$ and the equivariant
$\tau$ maps can be calculated using the equivariant Chern character
map, and the equivariant Todd class of the tangent bundle.  If
if $\beta \in K_0(Z,X^h)$ then using the formulas of Theorem 3.3 
\cite{EdGr:00}
\begin{equation} \label{eqn.rrr}
\tau^Z_{X^h}(\beta)  =  \ch^{Z}(\beta)\tau^Z([{\mathcal O_{X^h}}])
\end{equation}

If in addition, $Z$ is connected or $X^h$ has a $Z$-linearized ample
line bundle then
$$\tau^Z(\mathcal O_{X^h}) = \frac{\Td^Z(T_{X^h})}{\Td^Z(\z)}$$
where $\Td^Z$ is the equivariant Todd class of \cite[Definition 3.1]{EdGr:00}.
This follows from \cite[Theorem 3.1(d)]{EdGr:00} and
the following observation about the definition of
the equivariant Riemann-Roch map of \cite{EdGr:00}:
If a group $G$ acts freely on a smooth space $X$
with quotient 
$X \stackrel{\pi} \to X/G$
then, identifying $CH^*_G(X)$ with $CH^*(X/G)$, we have
$\tau^G({\mathcal O}_X) = \tau({\mathcal O}_{X/G})$
When $X/G$ (and thus $X$) is a smooth scheme then
$\tau({\mathcal O}_{X}) = \Td(T_{X/G})$. In this
case $\tau^G({\mathcal O}_X) = \Td^G(\pi^*T_{X/G})$. 
By Lemma \ref{lem.tantorsor}, $T_\pi = X \times {\mathfrak g}$.
Therefore, $\Td^G(\pi^*T_{X/G}) = \Td^G(T_X)/\Td^G({\mathfrak g})$.
\end{remark}

\begin{proof}
It suffices to prove Theorem \ref{t.rrtheorem1}, since this implies
Theorem \ref{t.rrtheorem2} using the Morita equivalence isomorphism
$G_0(G, S_{\Psi}) \simeq G_0(Z, X^h)$. We also need only
prove the second formula in Equation \eqref{e.rrtheorem1};
this implies the first formula because the maps $\tau$ and $\phi$
are covariant for finite morphisms.

The proof of Theorem
\ref{t.rrtheorem1} is almost the same as the proof in the case where
$G$ is diagonalizable, given in \cite[Theorem 3.1]{EdGr:03}.  The nonabelian
localization theorem replaces the localization theorem for actions of
diagonalizable groups used in \cite{EdGr:03}.  We give the general
proof here, referring to \cite{EdGr:03} for some omitted details.
Throughout the proof we write $\phi_M$ for $\phi^G_M$, where $G$ acts
properly on $M$.  We write $S$ for $S_{\Psi}$.

The proof proceeds in three steps.  First, we observe that the theorem
is true if the action of $G$ on $X$ is free (and that in this case, it
holds without the assumption that $X$ is smooth).  In this case, if
$\Psi \neq \{ 1 \}$, then $S$ is empty and $\alpha_{\Psi} = 0$, so both
sides of \eqref{e.rrtheorem1} vanish.  If $\Psi = \{1 \}$, then $f: S
\rightarrow X$ is an isomorphism, and the theorem amounts to the
assertion that
$$
\tau_Y(\alpha^G) = \phi_X \circ \tau^G_X(\alpha),
$$
which follows from \cite[Theorem 3.1(e)]{EdGr:00}.

Second, we prove the theorem for $\Psi = \{ 1 \}$.  In this case
$f\colon S \rightarrow X$ is an isomorphism, and the theorem amounts to the
assertion that
$$
\tau_Y((\alpha_1)^G) = \phi_X \circ \tau^G_X(\alpha).
$$
Since $\tau^G$ maps components of $K$-theory supported at maximal ideals
other than $\m_1$ to zero (this is proved as in \cite[Proposition~2.6]{EdGr:03}), 
on the right side of this equation we can replace $\alpha$ by $\alpha_1$.
By \cite[Proposition 10]{EdGr:98} there exists a finite surjective
morphism $p: X' \rightarrow X$ of $G$-spaces such that $G$ acts freely
on $X'$ and then induced map of quotients
$q\colon X' \to X$ is also finite and surjective. 
(Note that $X'$ need not be smooth.)  By \cite[Lemma 3.5]{EdGr:03}
(which holds without the diagonalizability
assumption on $G$), the map $p_*: G_0(G,X')_{\m_1} \rightarrow
G_0(G,X)_{\m_1}$ is surjective.  
Therefore there exists $\beta_1 \in
G_0(G,X')_{\m_1}$ such that $p_* (\beta_1) = \alpha_1$. There is a
geometric quotient $X' \rightarrow Y'$ such that the induced map $q:
Y' \rightarrow Y$ is also finite.  By Lemma \ref{lem.invcommute}, 
taking invariants commutes
with pushforward by a finite morphism, so $q_*((\beta_1)^G) = (\alpha_1)^G$.
Therefore,
$$
\tau_Y((\alpha_1)^G) = \tau_Y(q_*((\beta_1)^G) = q_* \tau_{Y'}((\beta_1)^G).
$$
By the first step of the proof, this equals $q_* \circ \phi_{X'} \circ
\tau^G_{X'}(\beta_1)$.  Since $q_* \circ \phi_{X'} = \phi_X \circ
p_*$, we have
$$
q_* \circ \phi_{X'} \circ \tau^G_{X'}(\beta_1) = \phi_X \circ p_*
\circ \tau^G_{X'}(\beta_1) = \phi_X \circ \tau^G_X \circ p_* (\beta_1)
= \phi_X \circ \tau^G_X (\alpha_1),
$$
as desired.

Third, we prove the theorem for general $\Psi$.  Let 
$\beta = \lambda_{-1}(N_f^*)^{-1} \cap f^* \alpha$.  
We need to prove that
$$
\tau_Y((\alpha_{\Psi}^G) = g_* \circ \phi_S \circ \tau^G_S (\beta (c_{\Psi})).
$$
By Lemma \ref{lemma.independenttwist}, 
the components of $\beta (c_{\Psi})$ and $\beta_{c_{\Psi}} (c_{\Psi})$ 
supported at $\m_1 \in R(G)$ are
equal, so arguing as in Step 2, we see that 
$\tau^G_S (\beta (c_{\Psi})) = \tau^G_S (\beta_{c_{\Psi}}(c_{\Psi}))$.
But $\beta_{c_{\Psi}} (c_{\Psi})$ is in $G_0(S,G)_{\m_1}$.  
Therefore we can apply
the second step of the proof to $\beta_{c_{\Psi}}(c_{\Psi})$.  Thus,
$$
\begin{array}{ccll}
\tau_Y(\alpha_{\Psi}^G) & = & 
\tau_Y \circ ((f_* (\beta_{c_{\Psi}}))^G) 
& \mbox{   localization} \\\\
 & = & g_* \circ \tau_{S/G} ( \beta_{c_{\Psi}}^G) 
& \mbox{  finite-pushforward commutes with invariants} \\\\
& = & g_* \circ \tau_{S/G} ( (\beta_{c_{\Psi}}(c_{\Psi}))^G) 
& \mbox{  Lemma \ref{l.ctwistinv}} \\\\
& = & g_*\circ \phi_S \circ \tau^G_S (\beta_{c_{\Psi}}(c_{\Psi}))
& \mbox { Step 2} \\\\
& = & g_* \circ \phi_S \circ \tau^G_S (\beta(c_{\Psi})). &
\end{array}
$$
This completes the proof.
\end{proof}

\section{Appendix}
This appendix contains a result about the tangent bundle to a
torsor which is difficult to find rigorously proved in the
literature.

\begin{lemma} \label{lem.tantorsor}
Let $X \stackrel{f} \to Y$ be a (left) $G$-torsor.
Then $T_f$ is canonically isomorphic to the $G$-bundle
$X \times \mathfrak{g}$ where the $G$-action
on the Lie algebra $\mathfrak{g}$ is the adjoint action.
\end{lemma}
\begin{proof}
By definition $T_f$ is the normal bundle to the diagonal
morphism $X \stackrel{\Delta_{f}} \to X \times_Y X$. 
Since $f$ is a $G$-torsor the diagram 
\begin{equation} \label{eq.diagtorsor}
\begin{array}{rrr}
G \times X & \stackrel{\sigma} \to &  X\\
\pi \downarrow & & f \downarrow\\
X  & \stackrel{f} \to & Y
\end{array}
\end{equation}
is a cartesian, where $\sigma \colon G \times X\to X$ is the action map
and $\pi\colon G \times X \to X$ is projection.  If $G$ acts on $G
\times X$ by conjugation on the first factor and the usual action on
the second factor then all morphisms in (\ref{eq.diagtorsor}) are
$G$-invariant. Thus there is a canonical identification of $G$-spaces
$X \times_Y X \to G \times X$. Under this identification the diagonal
corresponds to the section $X \stackrel{(e_G, 1_X)} \to G \times
X$. This map is obtained by base change from the $G$-equivariant inclusion
$e_G \to G$ (where $G$ acts on itself by conjugation)
whose normal bundle is $\mathfrak{g}$. Therefore, 
$T_f = X \times \mathfrak{g}$.
\end{proof}

\begin{remark}\label{rem.tantorsor}
Suppose $N \subset G$ is a closed normal subgroup with Lie algebra
$\mathfrak{n} \subset \mathfrak{g}$.
Normality of $N$ implies
that if $X \stackrel{f} \to Y$ is an $N$-torsor then there is
a natural left $G$-action on $Y$ such that $f$ is $G$-equivariant. 
Essentially
the same argument as above implies that $T_f$ is naturally
isomorphic to ${\mathfrak n} \otimes {\mathcal O}_X$ where
the $G$-action on ${\mathfrak n}$ is the restriction of
the adjoint action to the $ad$-invariant subalgebra ${\mathfrak n}\subset
{\mathfrak g}$.
\end{remark}

\def\cprime{$'$}


\begin{thebibliography}{MFK}

\bibitem[Ati]{Ati:74}
Michael~Francis Atiyah, {\em Elliptic operators and compact groups},
  Springer-Verlag, Berlin, 1974, Lecture Notes in Mathematics, Vol. 401.


\bibitem[BFM]{BFM:75}
Paul Baum, William Fulton, and Robert MacPherson, {\em Riemann-{R}och for
  singular varieties}, Inst. Hautes \'Etudes Sci. Publ. Math. (1975), no.~45,
  101--145.

\bibitem[Bor]{Bor:91}
Armand Borel, {\em Linear algebraic groups}, second ed., Springer-Verlag, New
  York, 1991.

\bibitem[BtD]{BrtD:95}
Theodor Br{\"o}cker and Tammo tom Dieck, {\em Representations of compact {L}ie
  groups}, Graduate Texts in Mathematics, vol.~98, Springer-Verlag, New York,
  1995, Translated from the German manuscript, Corrected reprint of the 1985
  translation.

\bibitem[CG]{ChGi:97}
Neil Chriss and Victor Ginzburg, {\em Representation theory and complex
  geometry}, Birkh\"auser Boston Inc., Boston, MA, 1997.


\bibitem[EG1]{EdGr:98}
Dan Edidin and William Graham, {\em Equivariant intersection theory}, Invent.
  Math. \textbf{131} (1998), no.~3, 595--634.

\bibitem[EG2]{EdGr:00}
\bysame, {\em Riemann-{R}och for equivariant {C}how groups}, Duke Math. J.
  \textbf{102} (2000), no.~3, 567--594.

\bibitem[EG3]{EdGr:03}
\bysame, {\em Riemann-{R}och for quotients and {T}odd classes of simplicial
  toric varieties}, Comm. in Alg. \textbf{31} (2003), 3735--3752.

\bibitem[Fog]{Fog:69}
John Fogarty, {\em Invariant theory}, W. A. Benjamin, Inc., New York-Amsterdam,
  1969.

\bibitem[Ful]{Ful:84}
William Fulton, {\em Intersection theory}, Springer-Verlag, Berlin, 1984.

\bibitem[EGA4]{EGA4}
A.~Grothendieck and J.~Dieudonn\'e, {\em \'{E}lements de {G}\'eom\'etrie
  {A}lg\'ebrique {I}{V}. \'{E}tude locale des schemas et des morphismes de
  sch\'emas}, Inst. Hautes \'Etudes Sci. Publ. Math. No. \textbf{20, 24, 28,
  32} (1964, 1965, 1966, 1967).

\bibitem[Hoc]{Hoc:65}
G.~Hochschild, {\em The structure of {L}ie groups}, Holden-Day Inc., San
  Francisco, 1965.

\bibitem[Hum]{Hum:95}
James~E. Humphreys, {\em Conjugacy classes in semisimple algebraic groups},
  Mathematical Surveys and Monographs, vol.~43, American Mathematical Society,
  Providence, RI, 1995.

\bibitem[KM]{KeMo:97}
Se{\'a}n Keel and Shigefumi Mori, {\em Quotients by groupoids}, Ann. of Math.
  (2) \textbf{145} (1997), no.~1, 193--213.

\bibitem[Knu]{Knu:71}
Donald Knutson, {\em Algebraic spaces}, Springer-Verlag, Berlin, 1971, Lecture
  Notes in Mathematics, Vol. 203.

\bibitem[MFK]{MFK:94}
D.~Mumford, J.~Fogarty, and F.~Kirwan, {\em Geometric invariant theory}, third
  ed., Springer-Verlag, Berlin, 1994.

\bibitem[Nie]{Nie:74}
H.~Andreas Nielsen, {\em Diagonalizably linearized coherent sheaves}, Bull.
  Soc. Math. France \textbf{102} (1974), 85--97.


\bibitem[OV]{OnVi:90}
A.~L. Onishchik and {\`E}.~B. Vinberg, {\em Lie groups and algebraic groups},
  Springer Series in Soviet Mathematics, Springer-Verlag, Berlin, 1990,
  Translated from the Russian and with a preface by D. A. Leites.

\bibitem[Qui]{Qui:73}
Daniel Quillen, {\em Higher algebraic {$K$}-theory. {I}}, Algebraic $K$-theory,
  I: Higher $K$-theories (Proc. Conf., Battelle Memorial Inst., Seattle, Wash.,
  1972), Springer, Berlin, 1973, pp.~85--147. Lecture Notes in Math., Vol. 341.

\bibitem[Ros]{Ros:02}
Wulf Rossmann, {\em Lie groups}, Oxford Graduate Texts in Mathematics, vol.~5,
  Oxford University Press, Oxford, 2002, An introduction through linear groups.

\bibitem[Seg1]{Seg:68a}
\bysame, {\em The representation ring of a compact {L}ie group}, Inst. Hautes
  \'Etudes Sci. Publ. Math. (1968), no.~34, 113--128.

\bibitem[Seg2]{Seg:68b}
Graeme Segal, {\em Equivariant {$K$}-theory}, Inst. Hautes \'Etudes Sci. Publ.
  Math. (1968), no.~34, 129--151.


\bibitem[Sri]{Sri:96}
V.~Srinivas, {\em Algebraic {$K$}-theory}, second ed., Progress in Mathematics,
  vol.~90, Birkh\"auser Boston Inc., Boston, MA, 1996.

\bibitem[Tho1]{Tho:86d}
R.~W. Thomason, {\em Comparison of equivariant algebraic and topological
  {$K$}-theory}, Duke Math. J. \textbf{53} (1986), no.~3, 795--825.

\bibitem[Tho2]{Tho:86i}
\bysame, {\em Lefschetz-{R}iemann-{R}och theorem and coherent trace formula},
  Invent. Math. \textbf{85} (1986), no.~3, 515--543.

\bibitem[Tho3]{Tho:87}
\bysame, {\em Algebraic {$K$}-theory of group scheme actions}, Algebraic
  topology and algebraic $K$-theory (Princeton, N.J., 1983), Ann. of Math.
  Stud., vol. 113, Princeton Univ. Press, Princeton, NJ, 1987, pp.~539--563.

\bibitem[Tho4]{Tho:88}
\bysame, {\em Equivariant algebraic vs.\ topological {$K$}-homology
  {A}tiyah-{S}egal-style}, Duke Math. J. \textbf{56} (1988), no.~3, 589--636.

\bibitem[Tho5]{Tho:92}
\bysame, {\em Une formule de {L}efschetz en {$K$}-th\'eorie \'equivariante
  alg\'ebrique}, Duke Math. J. \textbf{68} (1992), no.~3, 447--462.

\bibitem[Toe]{Toe:99}
B.~Toen, {\em Th\'eor\`emes de {R}iemann-{R}och pour les champs de
  {D}eligne-{M}umford}, $K$-Theory \textbf{18} (1999), no.~1, 33--76.

\bibitem[VV]{VeVi:02}
Gabriele Vezzosi and Angelo Vistoli, {\em Higher algebraic {$K$}-theory of
  group actions with finite stabilizers}, Duke Math. J. \textbf{113} (2002),
  no.~1, 1--55.

\end{thebibliography}
\end{document}